\documentclass[a4paper, 10pt, oneside]{amsart}

\usepackage[utf8]{inputenc}
\usepackage[T1]{fontenc}
\usepackage[english]{babel}
\usepackage{color}
\usepackage{amsmath} 
\usepackage{amssymb} 
\usepackage{amsthm}
\usepackage{graphicx}
\usepackage[colorlinks=true,linkcolor=blue]{hyperref}
\usepackage{cancel}

\usepackage{lmodern}
\usepackage{microtype}

\theoremstyle{plain}
\newtheorem{thm}{Theorem}[section]
\newtheorem*{thm*}{Theorem}
\newtheorem{prop}{Proposition}[section] 
\newtheorem{lem}{Lemma}[section] 
\newtheorem{cor}{Corollary}[section]
\newtheorem*{cor*}{Corollary}

\newtheorem{defi}{Definition}[section]

\newtheorem{rem}{Remark}

\newcommand {\R} {\mathbb{R}} \newcommand {\Z} {\mathbb{Z}}
\newcommand {\T} {\mathbb{T}} \newcommand {\N} {\mathbb{N}}

\newcommand {\p} {\partial}
\newcommand {\dt} {\partial_t}

\usepackage{tikz}
\usetikzlibrary{shapes,arrows,calc}

\newcommand{\bd}{\beta}

\newcommand{\Ell}{\mathcal{E}_{t}}
\newcommand{\Ello}{\left(\left(g(y)(\frac{\p_{y}}{k}-it)\right)^{2}+ \frac{g(y)}{kr(y)}(\frac{\p_{y}}{k}-it)-\frac{1}{r^{2}(y)}\right)}

\begin{document}

\date{\today}
\author{Christian Zillinger}
\address{Mathematisches Institut, Universität Bonn, 53115 Bonn, Germany}
\email{zill@math.uni-bonn.de}

\title{On circular flows: linear stability and damping}

\begin{abstract}
In this article we establish linear inviscid damping with optimal decay rates
around 2D Taylor-Couette flow and similar monotone flows in an annular domain  $B_{r_{2}}(0) \setminus B_{r_{1}}(0) \subset \mathbb{R}^{2}$.
Following recent results by Wei, Zhang and Zhao~\cite{Zhang2015inviscid}, we
establish stability in weighted norms, which allow for a singularity formation
at the boundary, and additional provide a description of the blow-up behavior. 
\end{abstract}

\maketitle
\tableofcontents

\section{Introduction}
\label{sec:introduction}

In this article we consider the 
linear stability and long-time asymptotic behavior of circular flows in an annular domain
$(x,y) \in B_{r_{2}}(0) \setminus B_{r_{1}}(0)$.
Such two-dimensional flows can for example be established
experimentally in rotating cylinders, where the rotation is
sufficiently slow as to not cause a (three-dimensional) Taylor-Couette instability.

In this setting, radial vorticities
\begin{align}
  \begin{split}
    \omega(x,y)&=\omega(r), \\
    v(x,y)&= \p_r \psi e_\theta=
    \begin{pmatrix}
      -y \\ x
    \end{pmatrix}
    \frac{\psi'(r)}{r}, \\
    \psi''(r)+\frac{1}{r}\psi'(r)&=\omega(r),
  \end{split}
\end{align}
are stationary solutions of the incompressible 2D Euler equations.

\begin{figure}[h]
  \centering
  \includegraphics[width=0.4\linewidth,page=1]{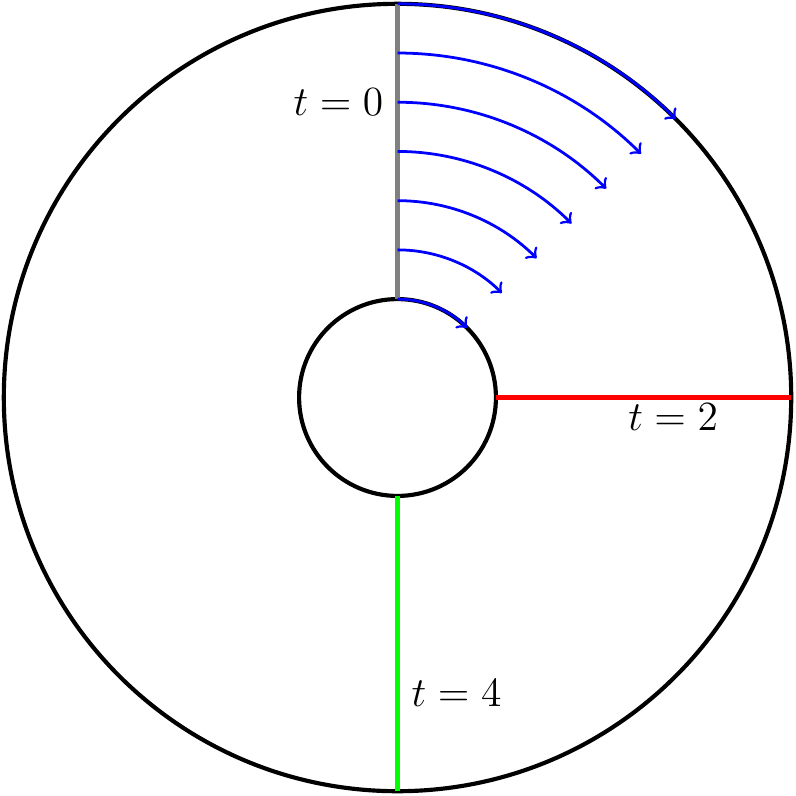}
  \includegraphics[width=0.5\linewidth, page=2]{pictures.pdf}
  \caption{Transport with constant angular velocity.
    We consider the Taylor-Couette flow $r$ in an annulus.
    The time $1$ flow-lines are drawn as arrows.
    A perturbation initially concentrated on a line stays concentrated on a
    line.
  On the right this behavior is expressed in polar coordinates.}
  \label{fig:TC}
\end{figure}

\begin{figure}[h]
  \centering
  \includegraphics[width=0.4\linewidth, page=3]{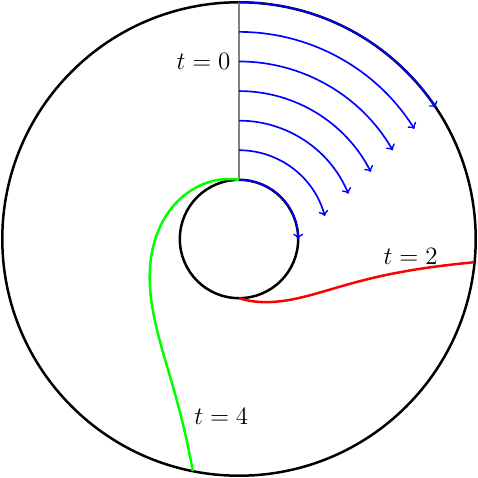}
  \includegraphics[width=0.4\linewidth, page=4]{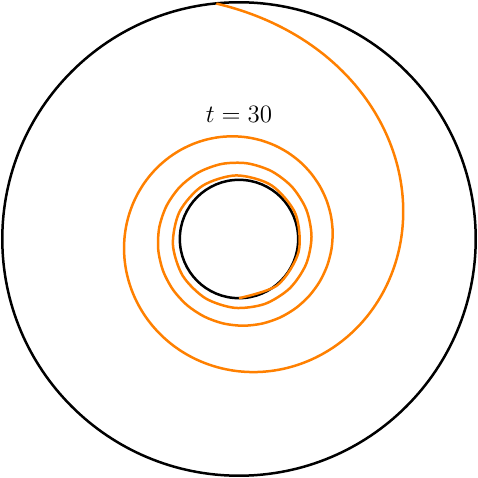}
  \includegraphics[width=0.4\linewidth, page=5]{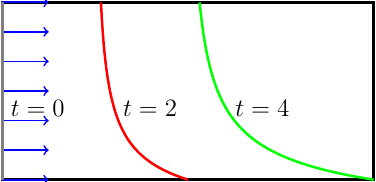}
  \includegraphics[width=0.4\linewidth, page=6]{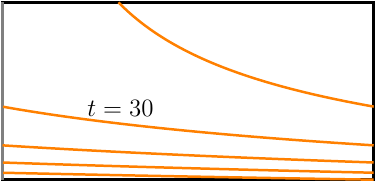}
  \caption{Transport by a \emph{monotone} flow.
    We consider the Taylor-Couette flow $r + \frac{1}{r}$, which we observe to
    be mixing. As time tends to infinity this mixing results in weak convergence to an averaged
    quantity.}
  \label{fig:Mon}
\end{figure}

Considering a small perturbation to Taylor-Couette flow,
\begin{align}
  \label{eq:1}
  \frac{\phi'(r)}{r}= A + \frac{B}{r^{2}},
\end{align}
we observe in Figure~\ref{fig:TC} that for $B=0$,
i.e. constant angular velocity, perturbations are rotated while
keeping their shape.  However, in the general case when $B \neq 0$,
$\frac{\phi'(r)}{r}$ is strictly monotone and the perturbation is sheared in way reminiscent of plane Couette flow,
as is depicted in Figure~\ref{fig:Mon}. This mixing behavior
underlies the phenomenon of (linear) inviscid damping.

Considering polar coordinates, the linearized Euler equations around
these stationary solutions are given by
\begin{align}
  \label{eq:polarLE}
  \begin{split}
    \dt f + U(r) \p_{\theta} f &= b(r) \p_{\theta} \phi, \\
    (\p_{r}^{2}+\frac{1}{r}\p_{r}+\frac{1}{r^{2}} \p_{\theta}^{2})\phi&= f, \\
    \p_{\theta}\phi|_{r=r_{1},r2_{2})}&=0, \\
    (t,\theta,r) & \in \R \times \T \times [r_{1},r_{2}],
  \end{split}
\end{align}
where $U$ and $b$ are given by
\begin{align*}
  U(r)&=\frac{\phi'(r)}{r}, \\
  b(r)&=-\frac{1}{r}\p_{r} (\p_{r}^{2}\phi(r)+\frac{1}{r}\p_{r}\phi(r)),
\end{align*}
and $b(r)\equiv 0$ if and only if one considers Taylor-Couette flow, $U(r)=A+\frac{B}{r^2}$.

As suggested by our notation, these equations share strong similarities with the linearized Euler
equations around a shear flow $(U(y),0)$
in a plane finite periodic channel, $\T \times [0,1]$:
\begin{align}
  \begin{split}
    \dt \omega + U(y) \p_{x} \omega - U''(y) \p_{x}\phi&=0, \\
    (\p_{y}^{2}+\p_{x}^{2})\phi &= \omega, \\
    \p_{x}\phi|_{y=0,1}&=0 , \\
    (t,x,y) &\in \R \times \T \times [0,1].
  \end{split}
\end{align}
Here, various different approaches have been used to study this and related settings.
\begin{itemize}
\item In \cite{stepin1995nonself}, Stepin studies the asymptotic stability of
  monotone shear flows using spectral methods. Under the assumption that the
  associated Rayleigh boundary value problem possesses no eigenvalues, he
  obtains an asymptotic description of the stream function and non-optimal decay
  rates.
\item In \cite{Euler_stability}, Bouchet and Morita provide heuristic results 
which suggest that the algebraic decay rates of Couette flow should hold for
general monotone flows as well.
However, their methods are not rigorous and do not provide sufficient error and
stability estimates, especially in higher Sobolev regularity, in order to prove
decay with optimal rates.
\item In \cite{Zill5} and \cite{Zill3}, the author establishes linear
  inviscid damping and scattering for monotone shear flows in an infinite and finite
  periodic channel. In the latter setting, we restrict to
  perturbations in $H^2 \cap H^{1}_0$ in order to obtain the optimal decay rates.
  Conversely, in the setting without vanishing Dirichlet boundary values, the
  sharp stability threshold is shown to be given by $H^{s},s=3/2$ due to
  asymptotic singularity formation at the boundary.
  
\item In \cite{Zhang2015inviscid}, Wei, Zhang and Zhao follow similar methods as in
  \cite{stepin1995nonself} and establish linear inviscid damping with optimal
  decay rates for monotone shear flows under the condition of there being no
  embedded eigenvalues.
  In particular, they remove the requirement of vanishing Dirichlet data and note that, using the boundary conditions of the velocity
  field and Hardy's inequality, one may allow for some blow-up at the boundary
  and still attain optimal decay rates.  
\item In a seminal work \cite{bedrossian2013asymptotic}, \cite{bedrossian2015inviscid} Bedrossian and Masmoudi establish nonlinear inviscid damping for Couette
  flow in an infinite periodic channel. There perturbations are required
  to be extremely regular, more precisely of Gevrey 2 class, in order to control
  nonlinear resonances.
  In particular, due to the singularity formation at the boundary and the
  associated blow-up of relatively low Sobolev norms, the question of linear
  inviscid damping for settings with boundary remains open.
\item In addition to the inviscid setting, Bedrossian,
  Germain and Masmoudi also consider Couette flow  as a solution of the
  Navier-Stokes equation in a two and three-dimensional infinite periodic
  channel.
  There, in addition to inviscid damping, the interaction between the mixing and
  viscous behavior yields additional stabilization by enhanced dissipation.
  Nonlinear inviscid damping is then established in Gevrey regularity \cite{bedrossian2015dynamics}
  and more recently in Sobolev regularity \cite{Bedrossian2015}, \cite{bedrossian2016sobolev}, where the
  threshold for stability results depends on $\nu>0$.
\item In the circular setting, research has focused on instability results,
  such as Taylor-Couette instability, bifurcation and turbulence. For an
  introduction we refer to the book of Chossat and Iooss
  \cite{chossat2012couette}. 
\end{itemize}

As the main results of this article we prove linear inviscid
damping and scattering for a general class of circular flows, satisfying
suitable monotonicity and smallness assumptions.
In comparison to our previous results, we note the following changes and improvements:
\begin{itemize}
\item We obtain optimal decay rates also for perturbations without vanishing
  Dirichlet data.
\item We show that $\p_yW$ splits into a
  bulk part $\Gamma$, which is stable also in unweighted higher Sobolev spaces,
and a boundary correction $\beta$, which is stable in a suitably weighted $H^1$
space, but exhibits blow-up in $L^{\infty}$.
\item The smallness condition is strongly reduced for results in higher
  regularity.
\item In this circular setting, periodicity in $\theta$ is a natural condition,
  unlike in the setting of a plane periodic channel.
\item We obtain a finer description of the boundary layer in terms of only the
  Dirichlet boundary values of the initial data.
\end{itemize}

\subsection{Main results}
\label{sec:main-results}
Our main results are summarized in the following theorem.

\begin{thm}[Linear inviscid damping with optimal decay rates]
  \label{thm:main1}
  Let $0<r_1<r_2<\infty$ and let $U:(r_1,r_2)\rightarrow (a,b)$ be bilipschitz and suppose that
  $h(\cdot)=b(U^{-1}(\cdot)) \in W^{3,\infty}((a,b))$ and that
  $\|h\|_{W^{1,\infty}}$ is sufficiently small.
  Then, for any $f_{0}\in H^{-1}_{\theta}H^{2}_{r}$ there exists $v_{\infty}(r)$
  such that the solution $f$ of~\eqref{eq:polarLE} satisfies 
  \begin{align}
   \label{eq:4}
    \|v(t,\theta,r)- v_{\infty}(r) e_{\theta}\|_{L^{2}} \lesssim <t>^{-1}\|f_{0}\|_{H^{-1}_{\theta}H^{1}_{r}}, \\
    \|v(t,\theta,r)e_{r}\|_{L^{2}} \lesssim <t>^{-2} \|f_{0}\|_{H^{-1}_{\theta}H^{2}_{r}},
   \end{align}
   as $t \rightarrow \infty$. There exists $f_{\infty} \in L^{2}_{\theta}H^{1}_{r}$ such that
   \begin{align*}
 f(t,\theta-tU(r),r) \rightarrow f_{\infty} \text{ in } L^{2}, 
   \end{align*}
and
    \begin{align}
      \label{eq:5}
    \|f(t,\theta-tU(r), r) - f_{\infty}(\theta, r)\|_{L^{2}_{\theta,r}} \lesssim <t>^{-1}\|f_{0}\|_{H^{-1}_{\theta}H^{2}_{r}}.
  \end{align}
\\

Furthermore, $f$ satisfies 
  \begin{align*}
    \|f(t,\theta-rU(r),r)\|_{H^{-1}H^{1}}+ \|(r-r_{1})(r-r_{2})\frac{d^{2}}{dr^{2}}f(t,\theta-rU(r),r)\|_{H^{-1}H^{1}} \lesssim \|f_{0}\|_{H^{-1}H^{2}}.
  \end{align*}
However, unless $bf|_{r=r_1,r_2}$ is constant,
\begin{align*}
  \sup_{t \geq 0} \|f(t,\theta-rU(r),r)\|_{H^{-1}H^{s}}= \infty,
\end{align*}
for any $s>3/2$.
More precisely, there exists an (explicit) function $\nu(t,\theta,r)$ determined solely by
$f_0|_{r=r_1,r_2}$ and $U$ such that 
\begin{align*}
  \|\frac{d^2}{dr^2} f(t,\theta-rU(r),r)- \nu \|_{L^2L^2} \leq \|f_0\|_{L^2H^2},
\end{align*}
and such that 
\begin{align*}
  \|(r-r_1)(r-r_2)\nu\|_{L^2} \leq |f_0|_{r=r_1,r_2}|.
\end{align*}
\end{thm}

\begin{rem}
  \begin{itemize}
  \item While $h=b(U^{-1})$ is required to be regular, the smallness assumption is only imposed on the $W^{1,\infty}$ norm. 
  \item This theorem summarizes the main results of Proposition
   ~\ref{prop:damping} and Theorems~\ref{thm:L2},~\ref{thm:StabilityofGamma} and
   ~\ref{thm:Stabilityofbeta} in terms of common norms in the variables
    $t,\theta,r$.
    In Section~\ref{sec:chang-vari-auxil} we introduce a \emph{scattering
    formulation}, which is used throughout the article. 
  \item The function $\nu$ is introduced in Section~\ref{sec:it-works}.
  \item In~\cite{Zhang2015inviscid} it has been observed that, by a use of
    Hardy's inequality, the second derivative of $W$ can be allowed to form a
    singularity as $t \rightarrow \infty$ while still attaining the optimal $t^{-2}$ decay rate.
    Here, we stress that stability in $H^2$ indeed does not hold due to singularity formation at the
    boundary as $t\rightarrow \infty$, as quantified in $\nu$ and $\beta$ (c.f. Section \ref{sec:Higher
      stability}).
  \item As we discuss in Section~\ref{sec:chang-vari-auxil}, our method of proof
    does not rely on cancellations or conserved quantities. Hence, the results extend to
    complex-valued $b(r)$ and various modified equations in a
    straightforward manner. 
    In the case of the linearized Euler equations in a plane finite
    periodic channel, however, Wei, Zhang and Zhao~\cite{Zhang2015inviscid} have
    shown, using different methods, that weaker assumptions suffice to obtain damping.
  \end{itemize}
\end{rem}

Similarly to~\cite{Zill3} our strategy is to first establish the damping and
scattering result, \emph{assuming} stability in higher Sobolev norms.
We stress that the damping estimate necessarily loses regularity.
Hence, usual Duhamel fixed point iteration approaches or energy
methods can not yield stability results.  Instead we employ a finer
study of the damping mechanism, which allows us to construct a
Lyapunov functional using the mode-wise decay to avoid the necessary
loss of regularity of uniform damping estimates.

The remainder of the article is organized as follows:
  \begin{itemize}
  \item In Section~\ref{sec:damping-mixing}, we show that regularity of the vorticity in
    coordinates moving with the flow can be exchanged for uniform damping
    estimates and that the problem of linear inviscid damping thus reduces to a
    stability problem.
    As motivating examples, we discuss the specific cases of Taylor-Couette flow, a
    point vortex and of Couette flow in a plane channel, where explicit
    solutions are available and, in a sense, trivial. 
  \item In Section~\ref{sec:chang-vari-auxil}, we introduce several reductions
    and changes of variables to arrive at a \emph{scattering formulation} of the
    linearized Euler equations.
    Subsequently, we analyze the structure of the equation and establish $L^{2}$ stability.
  \item Section~\ref{sec:Higher stability} considers higher regularity
    and singularity formation at the boundary. Compared to~\cite{Zill4}, in
    addition to considering a circular setting, we introduce a splitting  $\p_{y}W=\Gamma+\beta$, where $\Gamma$ is shown to be
    stable in higher regularity, regardless of Dirichlet boundary data.
    On the other hand, $\beta$ is determined solely by the underlying circular flow and the
    Dirichlet boundary data of the initial perturbation and
    provides an explicit characterization of the boundary layer.
    Subsequently, we further split $\p_y\beta$ to obtain an explicit
    characterization of the $H^2$ blow-up in the form of $\nu$ and stability in
    weighted spaces.
    Here, we rely on a new approach based on Duhamel's principle and an
    iterative estimate in order to control the evolution of the weighted quantities.
  \item The Appendices~\ref{sec:auxil-funct-bound} and~\ref{sec:Duhamel}
    provide a description of boundary evaluations for elliptic ODEs and a variant of
    Duhamel's formula adapted to a time-dependent right-hand-side of the equation.
  \end{itemize}

  \section{Damping by mixing, the role of regularity and examples}
  \label{sec:damping-mixing}

  As in the case of inviscid damping in a plane channel or
  Landau damping, decay of the velocity/force field and regularity of the
  solution in a coordinate system moving with the flow are closely linked.
  More precisely, in this section we show that uniform damping estimates closely
  correspond to a
  control of the regularity of 
  \begin{align*}
    W(t,\theta,r):=f(t,\theta-tU(r),r)
  \end{align*}
  with respect to $r$ and that such a control is necessary.
  The problem of linear inviscid damping with optimal decay rates thus turns out
  to be a stability problem, studied in Section \ref{sec:Higher stability}, which is the main focus of this article.
  \\
  
  We consider the linearized Euler equations
  \begin{align}
    \label{eq:polarLE_3}
    \begin{split}
      \dt f + U(r) \p_{\theta} f &= b(r) \p_{\theta} \phi, \\
      (\p_{r}^{2}+\frac{1}{r}\p_{r}+\frac{1}{r^{2}} \p_{\theta}^{2})\phi&= f, \\
      \p_{\theta}\phi|_{r=r_{1},r_{2}}&=0, \\
      (t,\theta,r) & \in \R \times \T \times [r_{1},r_{2}],
    \end{split}
  \end{align}
  as a perturbation around the transport problem
  \begin{align}
    \label{eq:transport_r}
    \begin{split}
    \dt g + U(r) \p_{\theta} g &=0, \\
    (t,\theta,r) & \in \R \times \T \times [r_{1},r_{2}].
    \end{split}
  \end{align}
  Based on this view, we measure the deviation of these equations by 
  introducing the \emph{scattered vorticity}
  \begin{align}
    W(t,\theta,r):= f(t,\theta-tU(r),r).
  \end{align}
  \emph{Assuming} regularity of $W$ uniformly in time, damping results for
  \eqref{eq:polarLE_3} then reduce to estimates for \eqref{eq:transport_r}.
  Here, we it has recently been observed by Wei, Zhang and Zhao
 ~\cite{Zhang2015inviscid} that quadratic decay rates only require
  control of a weighted $H^{2}$ norm
  \begin{align}
    \label{eq:9}
    \|W\|_{H^{1}}+ \|y (1-y)\p_{y}^{2}W\|_{L^{2}}
  \end{align}
  by using a Hardy inequality in the duality estimate.

  The following two propositions provide damping estimates in terms of
  regularity of $W$ in the case of a plane channel and a circular domain, respectively. 
  
  \begin{prop}[Damping by regularity for plane channel \cite{Zhang2015inviscid},
    \cite{Zill4}, \cite{Lin-Zeng}]
   \label{prop:damping_plane}
   Let $-\infty\leq a<b \leq \infty$ and let $U:(a,b) \rightarrow \R$ be
   locally $C^1$ and suppose that $U'(y)\neq 0$ for almost every $y\in (a,b)$.
   Let $W\in H^{-1}_xH^{1}_y(\T \times (a,b))$ with $\int_{\T}W dx=0$ and let 
   \begin{align*}
     \omega(t,x,y)=W(t,x-tU(y),y).
   \end{align*}
   Let further the associated velocity field $v$ be defined by
   \begin{align*}
     v_1&=-\p_y\phi,\\
     v_2&=\p_x\phi,\\
     \Delta \phi &=\omega, \\
     \p_x\phi|_{y=a,b}&=0, \\
     \phi &\in \dot H^1.
   \end{align*}
   Then $v$ satisfies 
   \begin{align}
     \label{eq:t1plane}
     \|v(t)\|_{L^2} &\lesssim \min \Big( \left\|W\right\|_{H^{-1}_{x}L^{2}_{y}}, t^{-1}\left\|\frac{W}{U'}\right\|_{H^{-1}_{x}H^{1}_{y}}, \\
& \quad t^{-1} \left(\left\|W|\frac{1}{U'}| + W|\p_y\frac{1}{U'}|\right\|_{H^{-1}_{x}L^{2}_{y}}+ \left\|(y-a)(y-b)\frac{\p_yW}{U'}\right\|_{H^{-1}_{x}L^{2}_{y}}\right)\Big).
   \end{align}
   Furthermore, suppose that $\p_{y}^2 W$ exists. Then, $v_2$ additionally satisfies 
   \begin{align}
     \label{eq:t2plane}
     \begin{split}
     \|v_2(t)\|_{L^2} &\leq t^{-2} \Big(\left\| \frac{W}{(U')^2}\right\|_{H^{-1}_{x}H^{1}_{y}} +  \left\| W\p_y\frac{1}{(U')^2}\right\|_{H^{-1}_{x}H^{1}_{y}}\\ & \quad + \min\left(\left\|\frac{(y-a)(y-b)}{(U')^2}\p_{y}^2W\right\|_{H^{-1}_{x}L^{2}_{y}},\left\|\frac{\p_{y}^2W}{(U')^2}\right\|_{H^{-1}_{x}L^{2}_{y}}\right)\Big).
     \end{split}
   \end{align}
  \end{prop}

  \begin{proof}[Proof of Proposition \ref{prop:damping_plane}]
  We note that, by integration by parts, 
  \begin{align*}
    \|v\|_{L^2}^2=\|\nabla \phi\|_{L^2}^2= -\iint \phi \omega dx dy.
  \end{align*}
  Applying Plancherel's theorem with respect to $x$ and noting that,
  \begin{align*}
    \mathcal{F}_x (\omega(t,\cdot,y))(k)= e^{iktU(y)} \hat{W}(t,k,y), 
  \end{align*}
  this equals 
  \begin{align*}
    \sum_{k \neq 0}\int \overline{\hat{\phi}}(t,k,y) e^{iktU(y)} \hat{W}(t,k,y). 
  \end{align*}
  Integrating 
  \begin{align*}
    e^{iktU}= \frac{1}{iktU'}\p_y e^{iktU(y)}
  \end{align*}
  by parts, we further obtain 
  \begin{align*}
    \sum_{k \neq 0} \frac{1}{t}\int e^{iktU(y)}\p_y \left(\overline{\hat{\phi}} \frac{\hat{W}}{ikU'}\right),
  \end{align*}
  which is controlled by 
  \begin{align*}
    t^{-1} \|\phi\|_{L^2H^1} \left\|\frac{W}{U'}\right\|_{H^{-1}H^1}.
  \end{align*}
  The estimate \eqref{eq:t1plane} thus follows by noting that 
  \begin{align*}
    \|\phi\|_{L^2H^1} \leq \|v\|_{L^2}.
  \end{align*}
\\

  In order to prove \eqref{eq:t2plane}, we note that 
  \begin{align*}
    \Delta v_2= \p_x \omega
  \end{align*}
  and define $\psi$ s.t.
  \begin{align*}
    \Delta \psi&=v_2, \\
    \p_x\psi|_{y=a,b}&=0, \\
    \psi &\in \dot H^{1}.
  \end{align*}
  Then, using integration by parts, we obtain 
  \begin{align*}
    \|v\|_{L^2}^2&= \iint \psi \p_x \omega= \sum_{k \neq} \overline{\hat{\psi}} e^{iktU(y)} ik\hat{W} dy \\
&= \frac{1}{t^2} \sum_{k \neq 0} \int e^{iktU(y)}\p_y \left(  \frac{1}{U'} \p_y \left(\frac{1}{U'} \overline{\hat{\psi}} \frac{\hat{W}}{k}\right)\right) dy \\
& \quad + \frac{1}{t^2} \sum_{k \neq 0} e^{iktU(y)}\frac{1}{U'} \p_y \left(\frac{1}{U'} \overline{\hat{\psi}} \frac{\hat{W}}{k}\right) \Big|_{y=a}^b \\
  \end{align*}
  The result hence follows by the Cauchy-Schwarz inequality, the trace map and by using the
  estimates
  \begin{align*}
    \|\phi\|_{H^1} &\lesssim \|v_2\|_{L^2}, \\
    \|\frac{\phi}{(y-a)(y-b)}\|_{L^2} &\lesssim \|\phi\|_{H^1},
  \end{align*}
  The first estimate here follows by standard elliptic regularity theory, while
  the second one is given Hardy's
  inequality, as observed in~\cite{Zhang2015inviscid}.
  \end{proof}

  The following proposition adapts these results to the setting of circular
  flows.
  \begin{prop}[Damping for circular flows;~\cite{Zhang2015inviscid},~\cite{Zill4},~\cite{Lin-Zeng}]
    \label{prop:damping}
    Let $0<r_1<r_2<\infty$ and let $U:(r_1,r_2)\rightarrow \R$ be locally $C^1$
    with $U'(r)\neq 0$ for almost every $r \in (r_1,r_2)$.  

    Let $\frac{W(t,\theta,r)}{U'} \in H^{-1}H^1(\T \times (r_1,r_2), r dr
    d\theta)$  with $\int_\T W d\theta=0$ and let 
    \begin{align*}
      f(t,\theta,r)=W(t,\theta-tU(r),r).
    \end{align*}
    Let further the associated velocity field be by defined by
    \begin{align}
      \begin{split}
        v_{r}(t,\theta,r)&= \frac{1}{r}\p_{\theta} \phi(t,\theta,r), \\
        v_{\theta}(t,\theta,r)&= \p_{r} \phi(t,\theta,r), \\
        (\p_{r}^{2}+ \frac{1}{r}\p_{r}
        +\frac{1}{r^{2}}\p_{\theta}^{2})\phi &=f, \\
        \p_{\theta}\phi|_{r=r_{1},r_{2}}&=0, \\
        v &\in L^2(r dr d\theta).
      \end{split}
    \end{align}
    Then $v$ satisfies
    \begin{align}
      \label{eq:t1circular}
      \|v(t)\|_{L^{2}(r dr d\theta)} &\lesssim \min \Big( \|W\|_{L^{2}(rdr d\theta)}, 
t^{-1} \|\frac{W(t)}{U'}\|_{H^{-1}_{\theta}H^{1}_{r}( r dr d\theta)}, \\
 & \quad t^{-1} \left(\|W(t)|\frac{1}{U'}|+ W r|\p_r\frac{1}{U'}|\|_{H^{-1}_{\theta}L^{2}_{r}( r dr d\theta)} + \left\| \frac{(r-r_1)(r-r_2)\p_rW}{U'}\right\|_{H^{-1}L^2(r dr d\theta)}\right) \Big).
    \end{align}
    Furthermore, suppose that $\p_{r}^2W$ exists. Then, $v_r$ additionally satisfies 
    \begin{align}
      \|v_{r}(t)\|_{L^{2}(r dr \theta)} &\lesssim t^{-2} \Big(\left\| \frac{W}{(U')^2}\right\|_{H^{-1}H^1(r dr d\theta)} + \left\| W \p_r \frac{1}{(U')^2}\right\|_{H^{-1}H^1( rdr d\theta)}   \\
& \quad + \min \left( \left\| \frac{\p_{r}^2 W}{(U')^2}\right\|_{H^{-1}L^2( r dr d\theta)}, \left\| \frac{(r-r_1)(r-r_2)\p_{r}^2 W}{(U')^2}\right\|_{H^{-1}L^2( rdr d\theta)} \right)\Big). 
    \end{align}
  \end{prop}
  
  \begin{rem}
  We note that for any given $0<r_1<r_2<\infty$ we could replace $r dr$ by just
  $dr$ in the above estimates at the cost of a constant $C(r_1,r_2)$.
  In this way the result and its proof can be made more similar to the setting
  of a finite channel.
  However, the above formulation also allows us to pass to the limits $r_1
  \downarrow 0$ and $r_2 \uparrow \infty$.
  \end{rem}

 \begin{proof}[Proof of Proposition~\ref{prop:damping}] 
  In order to obtain a more tractable stream function formulation of Euler's
  equations, in this proof we consider conformal coordinates, i.e. 
  \begin{align*}
    (r,\theta)&=(e^s,\theta),\\
 s \in (\log(r_1),\log(r_2))&=:(s_1,s_2).
  \end{align*}
  With respect to these coordinates, the stream function $\psi$ and the velocity
  field are given by 
  \begin{align*}
    e^{-2s}(\p_{s}^2+ \p_{\theta}^2)\psi &= \omega, \\
    v_r&=e^{-s}\p_\theta \psi, \\
    v_\theta&= e^{-s}\p_s \psi.
  \end{align*}
  Furthermore, the kinetic energy satisfies 
  \begin{align*}
    \int |v_r|^2 r dr d\theta = \int e^{-2s} |\p_\theta\psi|^2 e^{2s}ds d\theta &= \int |\p_\theta \psi|^2 ds d\theta, \\
    \int |v_\theta|^2 r dr d\theta &= \int |\p_s \psi|^2 ds d\theta, \\
    \int |v|^2 dr d\theta &= - \int \psi \omega e^{2s} ds d\theta.
  \end{align*}
  Applying a Fourier transform in $\theta$ and using the definition of $W$, we hence
  obtain 
  \begin{align*}
    \int |\p_s \psi|^2 + |\p_\theta\psi|^2 ds d\theta = - \sum_{k} \int \overline{\hat{\psi}} e^{2s} e^{iktU(e^s)} \hat{W} ds. 
  \end{align*}
  Integrating
  \begin{align*}
    e^{iktU(e^s)} = e^{-s} \frac{1}{iktU'(e^s)}\p_s e^{iktU(e^s)} 
  \end{align*}
  by parts, we further compute 
  \begin{align*}
    \int |\p_s \psi|^2 + |\p_\theta\psi|^2 ds d\theta=
 \sum_{k \neq 0} \int e^{iktU(e^s)} \p_s \left( \frac{1}{iktU'(e^s)} \overline{\hat{\psi}} e^{s} \hat{W}\right) ds.
  \end{align*}
  In order to estimate this integral, we use various different tools:
  \begin{itemize}
  \item If $\p_s$ does not fall on $\psi$, we control 
    \begin{align*}
     \sum_{k \neq 0} \int |\hat{\psi} \frac{1}{k}X| ds \leq \|\psi\|_{L^2(ds d\theta)} \|X\|_{H^{-1}_{\theta}L^{2}_s(ds d\theta)}  
    \end{align*}
    and use Poincar\'e's inequality to further estimate 
    \begin{align*}
\|\psi\|_{L^2(ds d\theta)} \leq C \|\p_\theta\psi\|_{L^2}. 
    \end{align*} 
  \item Alternatively, instead of Poincar\'e's inequality, duality yields an
    estimate by 
    \begin{align*}
      \|\p_\theta \psi\|_{L^2(ds d\theta)} \|X\|_{H^{-2}_{\theta}L^{2}_s(ds d\theta)}.
    \end{align*}
  \item Since $\psi$ has zero boundary values, we can also use Hardy's
    inequality to control by 
    \begin{align*}
      & \quad \|\frac{\psi}{(e^s-e^{s_1})(e^s-e^{s_2})}\|_{L^2(ds d\theta)} \|(e^s-e^{s_1})(e^s-e^{s_2}) X\|_{H^{-1}L^2(ds d\theta)}
\\ &\leq  \|\p_s \psi\|_{L^2(ds d\theta)} \|(e^s-e^{s_1})(e^s-e^{s_2}) X\|_{H^{-1}L^2(ds d\theta)}.
    \end{align*}
  \end{itemize}
In the case of a fixed annulus $\T \times (r_1,r_2),$ $ 0<r_1<r_2<\infty,$ the
precise choice of estimate is not essential.
However, when considering a non-periodic setting, e.g. $\R \times [a,b]$, or a
point vortex, i.e. $r_1=0$, or initial data with singularities at the boundary,
all these estimates can yield improvements.
\\

  In order to obtain the quadratic decay estimate for
  $v_{r}=e^{-s}\p_\theta\psi$, we note that 
  \begin{align*}
    (\p_{s}^2+\p_{\theta}^2)\p_\theta\psi= e^{2s}\p_{\theta}\omega.
  \end{align*}
  Thus, we define a potential $\gamma$ by 
  \begin{align*}
    (\p_{s}^2 + \p_{\theta}^2)\gamma &= \p_\theta \psi, \\
    \gamma|_{s=s_1,s_2}&=0,\\
    \nabla \gamma &\in L^2,                     
  \end{align*}
  and compute 
  \begin{align*}
    \int |v_r|^2 r dr d\theta &= \int |\p_\theta \psi|^2 ds d\theta \\
&= \int \gamma (\p_{s}^2+\p_{\theta}^2)\p_\theta\psi \\
&= \int \gamma e^{2s} \p_\theta \omega ds d\theta = \sum_{k\neq 0} ik \int \overline{\hat{\gamma}} e^{2s}e^{iktU(e^s)} \hat{W} ds.
  \end{align*}
  The result hence follows by integrating $e^{iktU(e^s)}$ by parts twice and
  using the Dirichlet data of $\gamma$ and $\p_\theta \psi$, the trace
  inequality and a variant of Hardy's inequality.
  That is, since $\gamma$ has zero Dirichlet boundary values, 
  \begin{align*}
    \|\frac{\gamma}{(e^s-e^{s_1})(e^s-e^{s_2})}\|_{L^2(ds d\theta)} &\leq \|\frac{\gamma}{(s-s_1)(s-s_2)}\|_{L^2(ds d\theta)} \\ &\lesssim \|\p_s \gamma\|_{L^{2}(ds d\theta)}= \|v_r\|_{L^2(rdr d\theta)}.
  \end{align*}
\end{proof}

We stress that these uniform damping estimates necessarily lose
regularity, since the associated change of coordinates is a unitary operator.
Thus, the operator norm of $f \mapsto v$ considered as a mapping from $L^2$ to
$L^2$ does not improve in time.
 Hence, it is not possible to derive stability of \eqref{eq:polarLE_3} using a
 common Duhamel-type approach or a
fixed point mapping. Instead, in Sections \ref{sec:L2stability} and \ref{sec:Higher stability} we have to make use of finer properties of the
dynamics and the mode-wise decay of the principal symbol of the
evolution operator.
Before that, in the following we discuss some examples for which explicit
computations are possible.

\subsubsection{Taylor-Couette flow}
\label{sec:examples}

As an application of the damping results, we discuss some
exceptional cases for which $W$ can be trivially computed in terms of the initial datum.

\begin{cor}[Couette flow]
Let $U(y)=y$ on $\T \times [a,b]$ with $a,b \in [-\infty,\infty]$,
then the linearized Euler equations reduce to the free transport equations.
Furthermore, if $\omega_0 \in H^{-1}_xH^{2}_y$, then the associated velocity
field satisfies 
\begin{align*}
  \|v(t)-\langle v|_{t=0} \rangle_x\|_{L^2}&\leq t^{-1}\|\omega_0\|_{H^{-1}H^1}, \\
  \|v_2(t)\|_{L^2} &\leq t^{-2}\|\omega_0\|_{H^{-1}H^2}.
\end{align*}
\end{cor}

\begin{cor}[Taylor-Couette flow; Point vortex]
\label{cor:pointvortex}
  Let $A, B \in \R$ and let $0\leq r_1<r_2 \leq\infty$, then the linearized Euler equations around Taylor-Couette
  flow 
  \begin{align*}
    U(r)= (Ar + \frac{B}{r})e_{\theta},
  \end{align*}
  are given by 
  \begin{align*}
    \dt f + (A+\frac{B}{r^2})\p_\theta f &=0, \text{ on } (0,\infty)\times \T \times (r_1,r_2) \\
    f|_{t=0}&=f_0 \text{ on } \T \times (r_1,r_2).
  \end{align*}
  Furthermore, the associated velocity field $v$ satisfies 
  \begin{align*}
    \|v\|_{L^2(rdrd\theta)} \leq C t^{-1}B^{-1} \|f_0\|_{H^{-1}H^1((r^7+r^5)drd\theta)}, \\
    \|v_r\|_{L^2( rdr d\theta)} \leq C t^{-2}B^{-2} \|f_0\|_{H^{-1}H^2((r^7+r^5)drd\theta)}. 
  \end{align*}
 Here, the case the case $r_1=0$, $A=0, B \neq 0$ corresponds to a
  \emph{point vortex}.
\end{cor}

\begin{proof}[Proof of Corollary \ref{cor:pointvortex}]
  We note that $(A+\frac{B}{r^2})'=-B \frac{2}{r^3}$.
  Hence, by direct computation 
  \begin{align*}
    \|\frac{\omega_0}{U'}\|_{H^{-1}H^1(rdrd\theta)} + \|\frac{\omega_0}{U'}\|_{L^2(r^{-1} dr \theta)} \leq 2B^{-1} (\|\omega_0\|_{H^{-1}L^2((r^7+r^5) dr \theta)} + \|\p_r\omega_0\|_{H^{-1}L^2(r^7 dr \theta)}).
  \end{align*}
\end{proof}

\section{Scattering formulation and $L^2$ stability}
\label{sec:chang-vari-auxil}

As established in Section~\ref{sec:damping-mixing}, 
the core problem of (linear) inviscid damping consists of establishing
a control of higher Sobolev norms of the vorticity moving with the flow: 
\begin{align}
  \label{eq:12}
W(t,\theta,r):= f(t,\theta-tU(r),r).
\end{align}
Here, we largely follow a similar approach as in the plane setting
considered in~\cite{Zill4}. As key improvements we obtain a less restrictive
smallness condition and develop a
splitting of $\p_{r}W$ into a well-behaved and more regular part $\Gamma$ and a (relatively) explicit
boundary layer $\bd$. This then allows us to deduce damping with optimal decay
rates and a detailed stability in suitable weighted Sobolev spaces, such as the ones considered
in Proposition~\ref{prop:damping}.

In order simplify our analysis, in this section we introduce
several changes of variables as well as useful auxiliary
functions.

\subsection{Scattering formulation}
\label{sec:scatt-form}
Expressing the linearized Euler equations
\begin{align}
  \label{eq:polarLE_2}
  \begin{split}
    \dt f + U(r) \p_{\theta} f &= b(r) \p_{\theta} \phi, \\
    (\p_{r}^{2}+\frac{1}{r}\p_{r}+\frac{1}{r^{2}} \p_{\theta}^{2})\phi&= f, \\
    \p_{\theta}\phi|_{r=r_{1},r_{2})}&=0, \\
    (t,\theta,r) & \in \R \times \T \times [r_{1},r_{2}],
  \end{split}
\end{align}
in terms of the \emph{scattered} quantities
\begin{align}
  \begin{split}
  F(t,\theta,r)&= f(t,\theta-tU(r),r), \\
  \Upsilon(t,\theta,r)&= \phi(t,\theta-tU(r),r),
  \end{split}
\end{align}
we obtain
\begin{align}
  \begin{split}
    \dt F &= b(r) \p_{\theta} \Upsilon, \\
    ((\p_{r}-tU'(r)\p_{\theta})^{2}+\frac{1}{r}(\p_{r}-tU'(r)\p_{\theta})+\frac{1}{r^{2}} \p_{\theta}^{2})\Upsilon&= F, \\
    \p_{\theta}\Upsilon|_{r=r_{1},r2_{2})}&=0, \\
    (t,\theta,r) & \in \R \times \T \times [r_{1},r_{2}],
  \end{split}
\end{align}
As none of the coefficient functions depend on $\theta$,
our system decouples with respect to Fourier modes $k$ in $\theta$.
\begin{align}
  \begin{split}
    \dt \hat{F} &= b(r) ik \hat{\Upsilon}, \\
    ((\p_{r}-iktU'(r))^{2}+\frac{1}{r}(\p_{r}-iktU'(r))-\frac{k^{2}}{r^{2}})\hat{\Upsilon}&= \hat{F}, \\
    ik\hat{\Upsilon}|_{r=r_{1},r2_{2})}&=0, \\
    (t,k,r) & \in \R \times 2\pi\Z \times [r_{1},r_{2}],
  \end{split}
\end{align}
We in particular note that the mode $k=0$,
which corresponds to a purely circular flow, is conserved in time.
Using the linearity of our equations, in the following we hence
without loss of regularity consider $k \in 2\pi (\Z\setminus \{0\})$
as a given parameter.

In view of the structure of the differential equation for $\Phi$,
it is further advantageous to use that $U$,
as a strictly monotone function, is invertible.  Introducing a change
of coordinates
\begin{align}
  r \mapsto y=U(r).
\end{align}
as well a denoting
\begin{align}
  \begin{split}
  h(y) &= \frac{(\omega_{0})'}{r}|_{r=U^{-1}(y)}, \\
  g(y)&= U'(r)|_{r=U^{-1}(y)}, \\
  W(t,y,k) &= \hat{F}(t,r,k)|_{r=U^{-1}(y)},\\
  \Phi(t,y,k) &= \frac{1}{k^{2}}\hat{\Upsilon}(t,r,k)|_{r=U^{-1}(y)},
  \end{split}
\end{align}
our system is then given by the following definition.
\begin{defi}[Euler's equations in scattering formulation]
  \label{defi:SLE}
  Let $U:[r_1,r_2]\rightarrow \R$ be strictly monotone and let $h(y)=b|_{r=U^{-1}(y)}$ and $g=U'(U^{-1}(y))$.
  Then \emph{Euler's equations in scattering
    formulation} are given by
  \begin{align}
    \label{eq:SLE}
    \begin{split}
      \dt W = \frac{ih(y)}{k}\Phi &=:  \frac{ih(y)}{k}L_{t}W, \\
      \Ell\Phi:= \Ello \Phi &=W, \\
      \Phi|_{y=a,b}&=0, \\
      (t,k,y) &\in \R \times 2\pi (\Z \setminus \{0\}) \times [a,b],
    \end{split}
  \end{align}
  where $a=\min(U^{-1}(r_{1}), U^{-1}(r_{2}))$, $b=\max(U^{-1}(r_{1}),
  U^{-1}(r_{2}))$ and $k \in 2 \pi (\Z \setminus \{0\})$.
\end{defi}

\begin{rem}
  \begin{itemize}
  \item Our methods do not rely on the specific form of $h$ or $g$ in terms of
    $U$. For example, we can allow for $h$ to be an arbitrary \emph{complex valued} $W^{1,\infty}$
    function.
  \item 
Here the notation $L_{t}W$
is used to stress that the mapping $W \mapsto \Phi$
is a linear operator in $W$.
\item 
As this system decouples with respect to $k$,
we will often treat $k\neq 0$
as a fixed given external parameter and with slight abuse of notation
use $W(t,y)$ to refer to $W(t,k,y)$ for the given $k$.
  \end{itemize}
\end{rem}

\subsection{Shifted elliptic regularity and modified spaces}
\label{sec:modif-spac-glid}

We note that in this scattering formulation $\Ell$ is obtained from an elliptic
operator by conjugation with $e^{ikty}$ and hence define suitable replacements
of the $H^1$ and $H^{-1}$ energies:

\begin{defi}[ $\tilde{H}^{1}_{t}$ and $\tilde{H}^{-1}_{t}$ energies]
  Let $u \in H^{1}([a,b])$ and let $k \in 2\pi(\Z \setminus \{0\})$ be given, then
  for every $t \in \R$, we define 
  \begin{align}
    \label{eq:tildeH1energy}
    \|u\|_{\tilde{H}^{1}_t}^2:= \|e^{ikty}u\|_{H^1}^2= \|u\|_{L^2}^2 + \|(\frac{\p_y}{k}-it)u\|_{L^2}^2. 
  \end{align}
  Furthermore, we define a dual quantity in the following way.
  Let $v \in L^2$ and let $\Psi[v]$ be the unique solution of 
  \begin{align*}
    (-1+(\frac{\p_y}{k}-it)^2) \Psi[v]&=v, \\
    \Psi[v]|_{y=a,b}&=0.
  \end{align*}
  Then we define 
  \begin{align*}
    \|v\|_{\tilde{H}^{-1}_t}:= \|\Psi[v]\|_{\tilde{H}^{1}_t}.
  \end{align*}
\end{defi}

\begin{lem}[Duality]
\label{lem:dual}
Let $W \in L^{2}$ and let $k \in (\Z\setminus\{0\})$ be given.
Then
  \begin{align}
    \label{eq:13}
    \|W\|_{\tilde{H}^{-1}_{t}}=\sup \{ \langle W, \alpha \rangle_{L^{2}}: \alpha \in H^{1}_{0}, \|\alpha\|_{\tilde{H}^{1}}\leq 1 \},
  \end{align}
i.e. $\tilde{H}^{-1}_{t}$ is dual to $\tilde{H}^{1}_{t}$. 
\end{lem}

\begin{proof}[Proof of Lemma~\ref{lem:dual}]
  Since multiplication by $e^{ikty}$ is a unitary operation and preserves zero
  Dirichlet boundary values and
  \begin{align*}
    \Psi_t[v]= e^{-ikty}\Psi_0[e^{ikty}v],
  \end{align*}
  it suffices to consider the case $t=0$, which is given by the usual $H^1$
  and $H^{-1}$ norms (where we use $\frac{\p_y}{k}$ instead of $\p_y$).

  The result then follows using integration by parts: 
  \begin{align}
    \label{eq:14}
    -\langle W, \alpha \rangle = \langle (1-\frac{\p_{y}}{k}^{2}) \Psi[W], \alpha \rangle \\
= \langle \Psi[W], \alpha  \rangle + \langle \frac{\p_{y}}{k}\Psi[W], \frac{\p_{y}}{k}\alpha \rangle
\leq \|W\|_{H^{-1}}\|\alpha\|_{H^{1}},
  \end{align}
  with equality if $\alpha= -\frac{1}{\|\Psi[W]\|_{H^1}}\Psi[W]$.
 Taking the supremum over all $\alpha$ with $\|\alpha\|_{H^{1}}$ we hence obtain
 the result.
\end{proof}

\subsection{Heuristics and obstructions}
\label{sec:heuristics}

On a \emph{heuristic} level, in order to establish stability in $L^2$, we use
that 
\begin{align*}
  \frac{d}{dt}\|W(t)\|_{L^2}^2 =2\Re \langle W, \frac{ih}{k}L_t W\rangle \lesssim C(h,k)\|W(t)\|_{\tilde{H}^{-1}_t}^2,
\end{align*}
and that for fixed functions $u \in L^2$, which \emph{do not depend on time}, 
\begin{align*}
  \int_{0}^{\infty} \|u\|_{\tilde{H}^{-1}_t}^2 dt \leq C \|u\|_{L^2}^2,
\end{align*}
as can be computed from a Fourier characterization.
Hence, it seems reasonable to expect that solutions $W(t)$ of~\eqref{eq:SLE}
satisfy an estimate of the form 
\begin{align*}
  \|W(t)\|_{L^2} \leq \exp(C\|h\|_{L^\infty}|k|^{-1}) \|f_0\|_{L^2},
\end{align*}
also for complex valued $h$, which is the case for some explicit model problems (c.f.~\cite{Zill3}). 

However, we stress that this heuristic is very rough and does not account for several obstructions:
\begin{itemize}
\item We note that integrability in time in general fails for time-dependent
$u \in L^{\infty}_t(L^2)$. For example, choosing
\begin{align*}
  u(t,k,y)=e^{ikty}u_0(k,y),
\end{align*}
we observe that 
\begin{align*}
  \int_{0}^{T} \|u\|_{\tilde{H}^{-1}_t}^2 dt = T \|u_0\|_{H^{-1}}^2,
\end{align*}
which diverges as $T \rightarrow \infty$ despite $\|u(t,y)\|_{L^2}= \|u_0\|_{L^2}$ being uniformly bounded.
\item Since the first estimate does not account for antisymmetric operators in
  $\frac{d}{dt}W$ it is not sufficient to establish $L^2$ stability.
For example, this estimate is satisfied by solutions $u(t,y)$ to 
\begin{align*}
  \dt u +iy u &= \Phi, \\
  (-1+(\p_y-it)^2)\Phi&=u, \\
  \Phi|_{y=a,b}&=0.
\end{align*}
Considering $v(t,y)=e^{ity}u(t,y)$, we observe that $v$ solves
\begin{align*}
  \dt v &= \phi, \\
  (1-\p_{y}^2)\phi &= v, \\
  \phi|_{y=a,b}&=0.
\end{align*}
Hence, choosing $u|_{t=0}$ to be an eigenfunction of $(1-\p_{y}^2)$, we obtain
an \emph{exponentially growing} solution.
\end{itemize}
 
\subsection{$L^2$ stability}
\label{sec:L2stability}

As the main result of this section, we adapt
the Lyapunov functional approach of~\cite{Zill5} to this circular
setting and prove stability of~\eqref{eq:SLE}.
In the following we formulate the main ingredients of our approach as a series
of Lemmata, which are then used to prove $L^2$ stability in Theorem
\ref{thm:L2}.
Subsequently, we elaborate on the theorem's statement and assumptions in comparison to
existing results and prove the lemmata.
Here, the lemmata are formulated in a general way in order to facilitate their
use for higher regularity estimates in later sections.

\begin{lem}
\label{lem:EllipticContribution}
  Let $L_{t}$ be given by~\eqref{eq:SLE} and let $\kappa \in W^{1,\infty}$.
  Then, for any $u,v \in L^{2}$ 
  \begin{align}
    \label{eq:15}
    | \langle u, \kappa L_{t}v \rangle | \leq (\|\kappa\|_{L^\infty} + \frac{1}{|k|} \|\p_y\kappa\|_{L^\infty}) \|u\|_{\tilde{H}^{-1}_{t}}\|L_{t} v\|_{\tilde{H}^{1}_t}
  \end{align}
\end{lem}

\begin{lem}
\label{lem:compCC}
  Let $L_{t}$ be as in~\eqref{eq:SLE}. Then there exists a constant $C=C(a,b,g)$ such that for any $u \in L^{2}$ and any $t\geq 0$
  \begin{align}
    \label{eq:16}
    \|L_{t}u\|_{\tilde{H}^{1}_{t}} \leq C \|u\|_{\tilde{H}^{-1}_{t}}.
  \end{align}
\end{lem}

\begin{lem}[{\cite[Lemma 4.5]{Zill4}}]
\label{lem:A}
  Let $u \in L^2([a,b])$ and let $\sum_{n \in (b-a)\N} u_{n} \sin(ny)$ be its series expansion.
  Define the symmetric, positive definite, non-increasing operator $A$ by 
  \begin{align}
    \label{eq:17}
    \langle u, A u \rangle := \sum_{n} \exp(\arctan(\frac{n}{k}-t)) |u_{n}|^{2}.
  \end{align}
  Then $A$ is symmetric, positive definite, non-increasing, $C^{1}$ in time and comparable to the identity, i.e.
  \begin{align}
    \label{eq:18}
    e^{-\pi} \|u\|_{L^{2}} \leq \langle u, A u \rangle \leq e^\pi \|u\|_{L^{2}},
  \end{align}
  for all $u \in L^{2}$.

  Furthermore, there exists a constant $e^{-\pi}\leq C_{2} \leq e^{\pi}$ and $\delta>0$ such that 
\begin{align}
  \label{eq:19}
  \|u\|_{\tilde{H}^{-1}_{t}}\|Au\|_{\tilde{H}^{-1}_{t}}  \leq -C_{2} \langle u, \dot A u \rangle 
\end{align}
\end{lem}

Using the preceding lemmata, we can establish $L^{2}$ stability.

\begin{thm}[$L^{2}$ stability]
\label{thm:L2}
Let $A$ and $C_{2}$ be given by Lemma~\ref{lem:A} and let $C$ be as in Lemma~\ref{lem:compCC}.
Further suppose that there exists $\delta>0$ such that
\begin{align}
  \label{eq:20}
  |k|^{-1}(\|h\|_{L^{\infty}}+|k|^{-1}\|\p_y h\|_{L^\infty}) \leq \frac{1}{C}(\frac{1}{C_2}-\delta).
\end{align}
Then for any solution $W$ to the Euler equations in scattering formulation~\eqref{eq:SLE} the functional
  \begin{align}
    \label{eq:21}
    I(t):= \langle W, A(t) W \rangle.
  \end{align}
  is non-increasing and satisfies 
  \begin{align}
    \label{eq:22}
    \dt I(t) \leq -\delta \|W(t)\|_{\tilde{H}^{-1}_t}\|A(t)W(t)\|_{\tilde{H}^{-1}_t} \leq 0.
  \end{align}
  In particular, this implies 
  \begin{align}
    \label{eq:23}
    e^{-\pi}\|W(t)\|_{L^{2}}^{2} \leq I(t) \leq I(0) \leq e^\pi \|\omega_{0}\|_{L^{2}}^{2}.
  \end{align}
\end{thm}

\begin{proof}[Proof of Theorem~\ref{thm:L2}]
  Using Lemma~\ref{lem:EllipticContribution}, we estimate
  \begin{align}
    \label{eq:24}
    \dt I(t) \leq \langle W, \dot A W \rangle + 2 |k|^{-1}(\|h\|_{L^{\infty}}+|k|^{-1}\|\p_y h\|_{L^\infty}) \|AW\|_{\tilde{H}^{-1}_{t}}\|L_tW\|_{\tilde{H}^{1}_{t}}.
  \end{align}
  Applying Lemma~\ref{lem:A} and Young's inequality, we further control
  \begin{align}
    \label{eq:25}
    \|AW\|_{\tilde{H}^{-1}_{t}}\|L_tW\|_{\tilde{H}^{-1}_{t}} \leq \frac{C}{2} \|AW\|_{\tilde{H}^{-1}_{t}} \|W\|_{\tilde{H}^{-1}_{t}}.
  \end{align}
  The result then follows by an application of Lemma~\ref{lem:A} and noting
  that, by our smallness assumption, 
  \begin{align*}
    \dt I(t) \leq \langle W, \dot A W \rangle + |k|^{-1}(\|h\|_{L^{\infty}}+|k|^{-1}\|\p_y h\|_{L^\infty})C\|AW\|_{\tilde{H}^{-1}_{t}}\|W\|_{\tilde{H}^{-1}_{t}} \\
\leq \langle W, \dot A W \rangle + (\frac{1}{C_2}-\delta)\|AW\|_{\tilde{H}^{-1}_{t}} \|W\|_{\tilde{H}^{-1}_{t}} \\
\leq -\delta \|AW\|_{\tilde{H}^{-1}_{t}}\|W\|_{\tilde{H}^{-1}_{t}} \leq 0
  \end{align*}
\end{proof}

Let us briefly remark on this result and its
assumptions: 
\begin{itemize}
\item We require a smallness condition on $\frac{ih}{k}L_t$ in order to rule out
the obstacles mentioned in Section~\ref{sec:heuristics}.
\item Since $h$ is allowed to be complex-valued, we do not rely on conserved
  quantities or classical stability results such as the ones of Rayleigh,
  Fjortoft or Arnold.
\item In the setting of a plane finite periodic channel,
  in~\cite{Zhang2015inviscid} Wei, Zhang and Zhao use a spectral approach to
  establish linear stability and decay with optimal rates for monotone shear flows under the assumption that the strictly monotone shear flow
  $U(y)$ possesses no embedding eigenvalues.
  In comparison, our smallness assumption is more restrictive, but extends to
  related problems such as stability in fractional Sobolev spaces, complex valued
  functions $h$ and fractional operators $L_t$ in a straightforward way.
\item In Section~\ref{sec:Higher stability}, we show that $\p_{y}^2W$ can be split into 
  a very regular, stable part $\Gamma$ and a boundary layer part $\beta$ which
  develops a singularity at the boundary.
  Here, $\beta$ is determined solely by the Dirichlet boundary data of the
  initial datum, $\omega_0$, and allows for a detailed study of the stability
  properties of the evolution.
\end{itemize}

It remains to prove Lemmata~\ref{lem:EllipticContribution} and~\ref{lem:compCC}.
\begin{proof}[Proof of Lemma~\ref{lem:EllipticContribution}]
  Let $u,v \in L^2$ and let $\Psi[u]$ be the unique solution of  
  \begin{align}
    \label{eq:26}
    (-1+(\p_{y}-ikt)^{2})\Psi[u]&=u, \\
    \Psi[u]|_{y=a,b}&=0.
  \end{align}
  Then we directly compute 
  \begin{align}
    \label{eq:27}
     | \langle u, \kappa L_{t}v \rangle | = | \langle  (1-(\frac{\p_{y}}{k}-ikt)^{2})\Psi[u], \kappa L_{t}v \rangle | \\
\leq \|\Psi[u]\|_{L^{2}} \|\kappa L_{t}\|_{L^{2}} + \|(\frac{\p_{y}}{k}-it)\Psi[u]\|_{L^{2}}\|(\frac{\p_{y}}{k}-it) \kappa L_{t}v\|_{L^{2}} \\
    \leq \|\Psi[u]\|_{\tilde{H}^{1}_{t}}(\|\kappa\|_{L^\infty} + \frac{1}{|k|} \|\p_y\kappa\|_{L^\infty})  \|L_{t}v\|_{\tilde{H}^{1}_{t}}.
  \end{align}
Here, we used that $L_tv$ by definition satisfies zero Dirichlet boundary
conditions and hence no boundary contributions appear when integrating by parts.
\end{proof}
\begin{proof}[Proof of Lemma~\ref{lem:compCC}]
  We recall that $L_{t}u$ is the solution of 
  \begin{align}
    \label{eq:28}
    \Ell L_{t} u= \Ello L_{t}u&=0 , \\
    L_{t}u|_{y=a,b}&=0,
  \end{align}
  and that $g(y)$ and $\frac{g(y)}{r(y)}$ are bounded from below (and above).

  Hence $\Ell$ is a shifted elliptic operator and testing by $L_{t}u$ (or $\frac{1}{g}L_{t}u$) we obtain that  
  \begin{align}
    \label{eq:29}
    \|L_{t}u\|_{\tilde{H}^{1}_{t}}^{2} \leq -C \langle L_{t}u, u \rangle,
  \end{align}
  for some $C>0$.
  Applying Lemma~\ref{lem:dual}, we thus obtain 
  \begin{align}
    \label{eq:30}
    \|L_{t}u\|_{\tilde{H}^{1}_{t}}^{2} &\leq C \|L_{t}u\|_{\tilde{H}^{1}_{t}}\|\Psi[u]\|_{\tilde{H}^{1}_{t}},\\
\Leftrightarrow \|L_{t}u\|_{\tilde{H}^{1}_{t}} &\leq C \|u\|_{\tilde{H}^{-1}_{t}}.
  \end{align}
\end{proof}

Having introduced the basic tools of our approach, in the following section we
consider higher stability of $W$, i.e. control of $\p_{y}W$.
Here, boundary effects qualitatively change the dynamics and necessitate a
modification of the weight $A(t)$.

\section{Higher stability and boundary layers}
\label{sec:Higher stability}

In this section we show that the $L^{2}$ stability result can be
extended to higher Sobolev regularity.  However, unlike in the setting
of an infinite periodic channel, boundary effects can not be neglected
and result in the formation of singularities.  As the main
improvements over our previous work for the plane channel in
\cite{Zill4}, we provide an explicit splitting into a more regular
good parts and a boundary layer exhibiting blow-up as well as an improved
smallness condition. 
This splitting then also allows to provide a more detailed description of the blow-up
also in weighted Sobolev spaces. For this purpose we
also introduce a different method of proof.

Let thus $W$ be a solution to \eqref{eq:SLE}
\begin{align*}
      \dt W &= \frac{ih}{k}L_{t}W, \\
      \Ell L_t W &=W, \\
      L_t|_{y=a,b}&=0, \\
      (t,k,y) &\in \R \times 2\pi (\Z \setminus \{0\}) \times [a,b].
\end{align*}
We begin by studying $\p_yW$, which satisfies 
\begin{align}
  \label{eq:pyW}
  \dt \p_yW &= \frac{ih}{k}L_t \p_y W + \frac{ih'}{k}L_tW + \frac{ih}{k}L_t[\Ell,\p_y]L_tW + \frac{ih}{k}H^{(1)}, \\
  \Ell H^{(1)}&=0 ,\\
  H^{(1)}_{y=a,b}&=\p_yL_tW|_{y=a,b}.
\end{align}
In contrast to the $L^2$ setting (or a setting without boundary such as $\T
\times \R$) we hence obtain a correction $H^{(1)}$ due to $\p_y L_tW$ not
satisfying zero Dirichlet boundary conditions.

As a main result of Appendix~\ref{sec:auxil-funct-bound}, we study the
boundary behavior of $\p_yL_t$ (also confer~\cite{Zill4}) and obtain the following description of $H^{(1)}$: 
\begin{lem}
  Let $W$ be a solution of~\eqref{eq:SLE} and let $H^{(1)}$ be the unique solution of 
  \begin{align*}
      \Ell H^{(1)}&=0 ,\\
  H^{(1)}_{y=a,b}&=\p_yL_tW|_{y=a,b}.
  \end{align*}
  Then there exist functions $u_1,u_2,\tilde{u}_1,\tilde{u}_2 \in H^{2}$
  (depending on $a,b,k$ and $g$ but not on $t$)
  and constants $c_1,c_2$ such that 
  \begin{align*}
    H^{(1)}(t,y)&= c_1\langle W, e^{ikt(y-a)}\tilde{u}_1 \rangle e^{ikt(y-a)}u_1 \\
    & \quad +  c_2\langle W, e^{ikt(y-b)}\tilde{u}_2 \rangle e^{ikt(y-b)}u_2. 
  \end{align*}
  Furthermore, for instance for $u_1$ for any $t>0$
  \begin{align*}
    \langle W, e^{ikt(y-a)}\tilde{u}_1 \rangle = \frac{\omega_0(a)}{ikt} - \frac{1}{ikt} \langle W, e^{ikt(y-a)}\p_y\tilde{u}_1 \rangle - \frac{1}{ikt} \langle \p_yW,  e^{ikt(y-a)}\tilde{u}_1  \rangle.
  \end{align*}
\end{lem}

Based on this characterization of $H^{(1)}$, we introduce a splitting of $\p_yW$
into a function $\beta$ depending only on $\omega_0|_{y=a,b}$ and
$\Gamma=\p_yW-\beta$.
As we show in Theorem~\ref{thm:StabilityofGamma}, $\Gamma$ is stable also in higher regularity.
In contrast, unless $\omega_0|_{y=a,b}$ is trivial, $\beta$ asymptotically
develops singularities at the boundary and exhibits blow-up in $H^{s},s>1/2$.
If one however considers weighted spaces, it is possible to compensate for these
singularities by vanishing weights and hence establish sufficient control for
damping with optimal decay rates.

\begin{lem}
\label{lem:splitting}
  Let $W$ be a solution of~\eqref{eq:SLE} and let $\Gamma$ be the
  solution of 
  \begin{align}
    \label{eq:Gamma}
    \begin{split}
    \dt \Gamma &= \frac{ih}{k}L_t \Gamma - \frac{h}{k^2t} \langle \Gamma,  e^{ikt(y-a)}\tilde{u}_1 \rangle e^{ikt(y-a)} u_1 - \frac{h}{k^2t} \langle \Gamma,  e^{ikt(y-b)}\tilde{u}_2 \rangle e^{ikt(y-b)} u_2 \\
 & \quad+ \frac{ih'}{k}L_tW + \frac{ih}{k}L_t[\Ell,\p_y]L_tW -c_1 \frac{h}{k^2t} \langle W, e^{ikt(y-a)}\p_y\tilde{u}_1 \rangle  e^{ikt(y-a)} u_1 \\ & \quad -c_2 \frac{h}{k^2t} \langle W, e^{ikt(y-b)}\p_y\tilde{u}_2 \rangle e^{ikt(y-b)} u_2, \\
\Gamma|_{t=0}&=\p_y \omega_0,
    \end{split}
  \end{align}
  and let $\beta$ be the solution of 
  \begin{align}
    \label{eq:beta}
    \begin{split}
    \dt \beta &= \frac{ih}{k}L_t \beta - \frac{h}{k^2t} \langle \beta,  e^{ikt(y-a)}\tilde{u}_1 \rangle e^{ikt(y-a)} u_1 - \frac{h}{k^2t} \langle \beta,  e^{ikt(y-b)}\tilde{u}_2 \rangle e^{ikt(y-b)} u_2 \\
    & \quad -\frac{c_1h\omega_0(a)}{k^2t} e^{ikt(y-a)}u_1  -\frac{c_2h\omega_0(b)}{k^2t} e^{ikt(y-b)}u_2, \\
    \beta|_{t=0}&=0.
    \end{split}
  \end{align}
  Then $\p_yW=\Gamma+\beta$.
  The function $\beta$ is called the \emph{boundary layer}.
\end{lem}

\begin{thm}[$H^2$ regularity of $\Gamma$]
  \label{thm:StabilityofGamma}
  Suppose that $g,h$ satisfy the assumptions of Theorem~\ref{thm:L2}.
  \begin{enumerate}
  \item 
    Suppose that additionally $g \in W^{2,\infty}$ and $h \in W^{2,\infty}$. 
  Then there exists a constant $C_1$ such that for all $\omega_0 \in H^1$ and any
  $t\geq 0$, the solution $\Gamma$ of~\eqref{eq:Gamma} satisfies
  \begin{align*}
    \|\Gamma(t)\|_{L^2} \leq C \|\omega_0\|_{H^1}.
  \end{align*}
\item 
  Suppose that additionally $g \in W^{3,\infty}$ and $h \in W^{3,\infty}$, then
  there exists a second constant $C_2$ such that for any 
  $\omega_0 \in H^2$ and for any $t \geq 0$, 
  \begin{align*}
        \|\Gamma(t)\|_{H^1} \leq C_2 \|\omega_0\|_{H^2}.
  \end{align*}
  \end{enumerate}
\end{thm}

\begin{thm}[$H^2$ regularity of $\beta$]
  \label{thm:Stabilityofbeta}
  Suppose $g,h$ satisfy the assumptions of Theorem~\ref{thm:L2}.
  \begin{enumerate}
  \item 
  Then there exists a constant $C_1$ such that for all $t\geq 0$, the solution
  $\beta$ of~\eqref{eq:beta} satisfies 
  \begin{align*}
    \|\beta(t)\|_{L^2} \leq C_1 (|\omega_0(a)|+|\omega_0(b)|).
  \end{align*}
  \item 
  Suppose that additionally $g,h \in W^{2,\infty}$, then there exists a second
  constant $C_2$ such that 
  \begin{align*}
     \|(y-a)(y-b)\p_y\beta(t)\|_{L^2} \leq C_2 (|\omega_0(a)|+|\omega_0(b)|).
  \end{align*}
  However, if for instance $|\omega_0(a)|>0$, then
  \begin{align*}
    |\beta(t,a)| \gtrsim \log(t)
  \end{align*}
  as $t \rightarrow \infty$ (similarly for $b$).
  In particular, by the Sobolev embedding, we obtain blow-up in $H^{s},s>1/2$.
  \end{enumerate}
\end{thm}

\begin{rem}
  \begin{itemize}
  \item  Combining Theorems~\ref{thm:L2},~\ref{thm:StabilityofGamma} and~\ref{thm:Stabilityofbeta}
    and Proposition~\ref{prop:damping}, we obtain Theorem~\ref{thm:main1}.
  \item It is possible to further split $\Gamma$ into functions controlled solely in
  terms of $\|\omega_0\|_{L^2}$, $\|\p_y \omega_0\|_{L^2}$ and $\|\p_{y}^2
  \omega_0\|_{L^2}$, if finer control is desired.  
\item Like Theorem~\ref{thm:L2}, in addition to these stability results we
  obtain Lyapunov functionals. As a key difference, these functionals are
  however in general only decreasing for times $t\geq T>0$. Control up time $T$
  is hence provided by a Gronwall-type argument, which determines the constants $C_1,C_2$. 
\item We stress that we do not require higher norms of $g,h$ to be small but
  only finite, so that
  derivatives of the equation are well-defined as mappings in $L^2$. 
\item When considering a setting without boundary contributions such as $\T
  \times \R$ or $\T \times \T$, no boundary correction $\beta$ is needed.
  Thus (a suitable modification of) this result already
  yields the desired stability for decay with optimal rates.
  Furthermore, this result generalizes to higher derivatives in a
  straightforward way, where again only finiteness of higher norms has to be required.  
  \end{itemize}
\end{rem}

\subsection{Stability of $\Gamma$}
\label{sec:h1gamma}
As the main result of this subsection we provide a proof of Theorem
\ref{thm:StabilityofGamma}.
Here, the $L^2$ stability result is self-contained, while the $H^1$ estimate
presupposes the $L^2$ stability of $\beta$, which is established in the
following subsection.
Furthermore, we briefly discuss the implications of Theorem
\ref{thm:StabilityofGamma} for settings
without boundary and provide an improved stability result for the setting of an
infinite plane periodic channel, $\T_L\times \R$.

We recall that $\Gamma$ is the solution of 
  \begin{align*}
    \dt \Gamma &= \frac{ih}{k}L_t \Gamma - \frac{h}{k^2t} \langle \Gamma,  e^{ikt(y-a)}\tilde{u}_1 \rangle e^{ikt(y-a)} u_1 - \frac{h}{k^2t} \langle \Gamma,  e^{ikt(y-b)}\tilde{u}_2 \rangle e^{ikt(y-b)} u_2 \\
 & \quad+ \frac{ih'}{k}L_tW + \frac{ih}{k}L_t[\Ell,\p_y]L_tW -c_1 \frac{h}{k^2t} \langle W, e^{ikt(y-a)}\p_y\tilde{u}_1 \rangle  e^{ikt(y-a)} u_1 \\ & \quad -c_2 \frac{h}{k^2t} \langle W, e^{ikt(y-b)}\p_y\tilde{u}_2 \rangle e^{ikt(y-b)} u_2, \\
\Gamma|_{t=0}&=\p_y \omega_0,
  \end{align*}
In addition to the estimates for $L_t$ derived in Section~\ref{sec:L2stability},
we hence need to control contributions of the form 
\begin{align*}
 \frac{1}{ikt} \langle \Gamma,  e^{ikt(y-a)}\tilde{u}_1 \rangle \langle A \Gamma, e^{ikt(y-b)} u_2 \rangle,
\end{align*}
which can not be controlled by the previous choice of $A(t)$.

Instead, we construct a modified weight $A_1(t)$, which is introduced in the
following Lemmata (cf. \cite{Zill5} for a similar construction adapted to
fractional Sobolev spaces).

\begin{lem}
\label{lem:modweight1}
  Let $u \in H^1$, then for $0<\mu<1/2$ and for
  every $v=\sum_n v_n e^{iny} \in L^2$ 
  \begin{align*}
    |\langle  v, e^{ikt(y-a)} u \rangle|^{2} \leq C_\mu \|u\|_{H^1}^2 \sum_n <n-kt>^{-2\mu}|v|_{n}^2.
  \end{align*}
\end{lem}

\begin{proof}[Proof of Lemma~\ref{lem:modweight1}]
  By expanding the $L^2$ inner product in a basis, we obtain that
  \begin{align*}
    \langle  v, e^{ikt(y-a)} u\rangle= \sum_n v_n \langle e^{iny},e^{ikt(y-a)} u \rangle.
  \end{align*}
  Integrating by parts and using the trace inequality, we further estimate 
  \begin{align*}
    |\langle  e^{iny} ,e^{ikt(y-a)} u \rangle | \leq <n-kt>^{-1} \|u\|_{H^1}.
  \end{align*}
  The result hence follows by an application of the Cauchy-Schwarz inequality:
  \begin{align*}
    |\langle  v, e^{ikt(y-a)} u\rangle | &\lesssim \sum_n v_n <n-kt>^{-\mu} <n-kt>^{1-\mu} \\
&\leq \|v_n <n-kt>^{-\mu}\|_{l^2} \|<n-kt>^{1-\mu}\|_{l^2} \\
&\leq 2 \|v_n <n-kt>^{-\mu}\|_{l^2} \|<n>^{-(1-\mu)}\|_{l^2} \\
&=: C_\mu \|v_n <n-kt>^{-\mu}\|_{l^2},
  \end{align*}
  where we used that $<n>^{-(1-\mu)} \in l^2$ if $\mu<1/2$.
\end{proof}

\begin{lem}
\label{lem:A1}
  Let $0<\lambda,\mu <1$ with $\lambda+2\mu>1$ and let $\epsilon>0$ and define the symmetric operator $A_1(t)$
  by its action on the basis:
  \begin{align*}
   A_1(t): e^{iny} \mapsto \exp\left(\arctan\left(\frac{\eta}{k}-t\right)- \epsilon \int^{t} <\tau>^{-\lambda} <n-k\tau>^{-2\mu} d\tau\right).
  \end{align*}
  Then for every $u \in L^2$ and every $t \in \R$ 
  \begin{align*}
    C \|u\|_{L^2}^2\leq \langle u, A_1(t) u \rangle &\leq C^{-1} \|u\|_{L^2}^2, \\
    \langle u, \dot{A}_1(t)u \rangle \leq -C_1\|u\|_{\tilde{H}^{-1}_t}^2 - C\epsilon \sum_n  <t>^{-\lambda} <n-kt>^{-2\mu}|u_n|^2 &\leq 0.
  \end{align*}
\end{lem}

\begin{proof}[Proof of Lemma \ref{lem:A1}]
  We note that $<t>^{-\lambda} <n-kt>^{-2\mu} \in L^1(\R)$ and that 
  \begin{align*}
    - \epsilon \int^{t} <\tau>^{-\lambda} <n-k\tau>^{-2\mu} d\tau
  \end{align*}
  is monotonically decreasing.
  The properties of $A_1(t)$ hence follow by direct computation, where 
  \begin{align*}
    C=\exp(-\pi - \epsilon \|<\cdot>^{-\lambda} <n-k\cdot>^{-2\mu}\|_{L^1(\R)}).
  \end{align*}
  and $C_1$ is determined by $C$ and Lemma~\ref{lem:A}.
\end{proof}

\begin{lem}
  \label{lem:com1}
  Let $g \in W^{2,\infty}$, $g \geq c>0$, then for every $u \in L^2$ and for
  every $t \geq 0$, 
  \begin{align*}
    \|L_t[\Ell, \p_y]L_tu\|_{\tilde{H}^{1}_{t}} \lesssim \|u\|_{\tilde{H}^{-1}_t}.
  \end{align*}
\end{lem}

\begin{proof}[Proof of Lemma \ref{lem:com1}]
  By Lemma~\ref{lem:compCC}, we obtain that 
  \begin{align*}
    \|L_t[\Ell, \p_y]L_tu\|_{\tilde{H}^{1}_{t}} \lesssim  \|[\Ell, \p_y]L_tu\|_{\tilde{H}^{-1}_{t}}.
  \end{align*}
  We further note that 
  \begin{align*}
    [\Ell, \p_y]= e^{-ikty}[\mathcal{E}_0, \p_y-ikt]e^{ikty}= e^{-ikty}[\mathcal{E}_0, \p_y]e^{ikty},
  \end{align*}
  and that, by direct computation, $[\mathcal{E}_0, \p_y]$ is a second-order
  operator.
  Hence, using integration by parts, we further estimate
  \begin{align*}
    \|[\Ell, \p_y]L_tu\|_{\tilde{H}^{-1}_{t}} \lesssim \|L_tu\|_{\tilde{H}^{1}_t} \lesssim \|u\|_{\tilde{H}^{-1}_t}.
  \end{align*}
\end{proof}

Using these results, we can now provide a proof of Theorem
\ref{thm:StabilityofGamma} and thus establish $L^2$ stability.

\begin{proof}[Proof of Theorem~\ref{thm:StabilityofGamma}, part  $(1)$]
Fix $0<\lambda,\mu<1$ with $2\mu+\lambda>1$ and let $A_1$ be given by Lemma
\ref{lem:A1}, where
\begin{align*}
0<\epsilon < \frac{1}{100} \|<n-k\cdot>^{-2\mu}<\cdot>^{-\lambda}\|_{L^1(\R)}^{-1}.
\end{align*}
Then we define 

\begin{align}
  \label{eq:LyapGamma1}
  I(t):= \langle \Gamma, A_1(t)\Gamma \rangle + C_1 \langle W, A(t)W \rangle,
\end{align}
where $C_1 \gg 0$ is to be chosen later.
We then claim that there exists $T>0$ such that for all initial data and for all
$t \geq 0$, $I(t)$ satisfies 
\begin{align*}
  \frac{d}{dt}I(t) \leq C t^{-2(1-\mu/2)} \|\omega_0\|_{L^2}^2 \in L^1(\R).
\end{align*}
Using Gronwall's inequality, we further obtain that 
\begin{align*}
  I(T)\leq \exp(CT)I(0),
\end{align*}
which concludes the proof.
\\

It remains to prove the claim. Using Theorem~\ref{thm:L2} and Lemma
\ref{lem:A1}, we directly compute 
\begin{align*}
  \frac{d}{dt}I(t) &\leq -C\|\Gamma\|_{\tilde{H}^{-1}_t}^2 - C\epsilon \sum_n  <t>^{-\lambda} <n-kt>^{-2\mu}|\Gamma_n|^2 \\ & \quad - C_1\delta \|W(t)\|_{\tilde{H}^{-1}_t}^2 + 2 \Re \langle \frac{d}{dt}\Gamma, A_{1}(t)\Gamma \rangle .
\end{align*}
Using Lemma~\ref{lem:modweight1}
and Lemma~\ref{lem:EllipticContribution} and recalling~\eqref{eq:Gamma}, we
further estimate 
\begin{align*}
  2 \Re \langle \frac{d}{dt}\Gamma, A_{1}(t)\Gamma \rangle &\leq C(h,k)\|\Gamma\|_{\tilde{H}^{-1}_t} \|A_1\Gamma\|_{\tilde{H}^{-1}_t} + C(h,k,\mu) \frac{1}{t} \left( \sum_n <n-kt>^{-2\mu}|\Gamma_n|^2 \right) \\
  & \quad + C(h,h',k) \|A_{1}\Gamma\|_{\tilde{H}^{-1}_t} (\|L_tW\|_{\tilde{H}^{1}_t} + \|L_t [\Ell,\p_y]L_t\|_{\tilde{H}^{1}_t}) \\
  & \quad + C(h,k,g) t^{-1} \|\omega_0\|_{L^2} \sqrt{\sum_n <n-kt>^{-2\mu}|\Gamma_n|^2}.
\end{align*}
Splitting $t=t^{-(1-\mu)} t^{-\mu}$ and using Young's inequality and Lemmata~\ref{lem:compCC}
and~\ref{lem:com1}, we further
control
\begin{align*}
  \frac{1}{t} \left( \sum_n <n-kt>^{-2\mu}|\Gamma_n|^2 \right)= t^{-(1-\lambda)}\sum_n t^{-\lambda} <n-kt>^{-2\mu}|\Gamma_n|^2 , \\
   \|A_{1}\Gamma\|_{\tilde{H}^{-1}_t} (\|L_tW\|_{\tilde{H}^{1}_t} + \|L_t [\Ell,\p_y]L_t\|_{\tilde{H}^{1}_t}) \leq \sigma \|A_{1}\Gamma\|_{\tilde{H}^{-1}_t}^{2} + \sigma^{-1} \|W(t)\|_{\tilde{H}^{-1}_t}^2 , \\
  t^{-1} \|\omega_0\|_{L^2} \sqrt{\sum_n <n-kt>^{-2\mu}|\Gamma_n|^2} \leq \sigma \langle \Gamma, \dot{A}_1(t) \Gamma \rangle| + \sigma^{-1} t^{-2(1-\mu/2)}\|\omega_0\|_{L^2}^2.
\end{align*}
Choosing $\sigma$ sufficiently small and letting $T>0$ be sufficiently large and using the smallness assumption of
Theorem~\ref{thm:L2}, we observe that  
\begin{align*}
  -C\|\Gamma\|_{\tilde{H}^{-1}_t}^2 - C\epsilon \sum_n  <t>^{-\lambda} <n-kt>^{-2\mu}|\Gamma_n|^2 
+ C(h,k)\|\Gamma\|_{\tilde{H}^{-1}_t} \|A_1\Gamma\|_{\tilde{H}^{-1}_t} \\
+ (C(h,k,\mu)t^{-(1-\mu)}) \sum_n <t>^{-\lambda} <n-kt>^{-2\mu}|\Gamma_n|^2 \\
+ \sigma \|A_{1}\Gamma\|_{\tilde{H}^{-1}_t}^{2}  + \sigma \langle \Gamma, \dot{A}_1(t) \Gamma \rangle| \leq 0.
\end{align*}
Similarly, choosing $C_1$ sufficiently large, we observe that 
\begin{align*}
  - C_1\delta \|W(t)\|_{\tilde{H}^{-1}_t}^2 + \sigma^{-1} \|W(t)\|_{\tilde{H}^{-1}_t}^2 \leq 0.
\end{align*}
Hence, we conclude that for $t \geq T>0$, $I(t)$ satisfies 
\begin{align*}
  \frac{d}{dt}I(t) \leq \sigma^{-1} t^{-2(1-\mu/2)}\|\omega_0\|_{L^2}^2.
\end{align*}
which finishes the proof of the claim and hence of the $L^2$ stability result, $(1)$.
\end{proof}

Next, we consider the
evolution of $\p_y\Gamma$:
\begin{align}
\label{eq:pyGamma}
  \dt \p_y\Gamma&= \frac{ih}{k}L_t \p_y \Gamma + \frac{ih'}{k}L_t \Gamma + \frac{ih}{k}L_t[\Ell,\p_y]L_t\Gamma - (\p_yL_t\Gamma)(a) e^{ikt(y-a)}u_1 - (\p_yL_t\Gamma)(b) e^{ikt(y-b)}u_1\\
& \quad+ \p_y\Big( \frac{h}{k^2t} \langle \Gamma,  e^{ikt(y-a)}\tilde{u}_1 \rangle e^{ikt(y-a)} u_1 - \frac{h}{k^2t} \langle \Gamma,  e^{ikt(y-b)}\tilde{u}_2 \rangle e^{ikt(y-b)} u_2 \\
 & \quad-c_1 \frac{h}{k^2t} \langle W, e^{ikt(y-a)}\p_y\tilde{u}_1 \rangle  e^{ikt(y-a)} u_1  -c_2 \frac{h}{k^2t} \langle W, e^{ikt(y-b)}\p_y\tilde{u}_2 \rangle e^{ikt(y-b)} u_2 \Big) \\
& \quad + \p_y \left( \frac{ih'}{k}L_tW + \frac{ih}{k}L_t[\Ell,\p_y]L_tW  \right), \\
\p_y \Gamma|_{t=0}&=\p_{y}^2 \omega_0
\end{align}
Since we here also have to compute $\p_yW=\Gamma+\beta$ in order to control $\|\p_y
\Gamma\|_{L^2}$, we require $L^2$ estimates on $\beta$.
Before continuing with the proof of Theorem~\ref{thm:StabilityofGamma}, we hence
prove the first part of Theorem~\ref{thm:Stabilityofbeta}
as well as some further properties of the evolution of $\beta$, which are formulated in the following proposition.

\begin{prop}
  \label{prop:betaL2}
  Suppose $g,h$ satisfy the assumptions of Theorem~\ref{thm:Stabilityofbeta}.
  Let $\beta$ be the solution of~\eqref{eq:beta} and let $A_1(t)$ be given by
  Lemma~\ref{lem:A1}. Then there exists $T>0$ such that for all $t\geq 0$
  \begin{align*}
    I_2(t)=\langle \beta, A_1(t)\beta \rangle
  \end{align*}
  satisfies 
  \begin{align*}
    \frac{d}{dt}I_2(t) \leq  \delta \langle \beta, \dot A_1(t)\beta \rangle + C t^{-2(1-\mu/2)}|\omega_0|_{y=a,b}|^2.
  \end{align*}
\end{prop}

\begin{proof}[Proof of Proposition~\ref{prop:betaL2}]
  Using the same weight $A_{1}$, we observe that 
  \begin{align*}
    \Re \langle A_{1}(t)\beta, \frac{ih}{k}L_t \beta -\frac{h}{k^2t} \langle \beta,  e^{ikt(y-a)}\tilde{u}_1 \rangle e^{ikt(y-a)} u_1 - \frac{h}{k^2t} \langle \beta,  e^{ikt(y-b)}\tilde{u}_2 \rangle e^{ikt(y-b)} u_2  \rangle \\
\leq (C(h,k)+ C(h,k,g)t^{-(1-\mu)})|\langle  \beta, \dot{A}_1(t)\beta \rangle|.
  \end{align*}
  Using the smallness assumption and restricting to $t\geq T>0$, this
  contribution can thus be absorbed by 
  \begin{align*}
    \langle  \beta, \dot{A}_1(t)\beta \rangle \leq 0.
  \end{align*}

  Hence, we focus on 
  \begin{align*}
    \Re \langle A_{1}\beta, \omega_{0}(a)\frac{1}{ikt}e^{ikty}u \rangle \lesssim C_{\lambda} \|\beta_{n}<n-kt>^{-\lambda}\|_{l^{2}} |\omega_{0}|_{y=a,b}| \frac{1}{|kt|}.
  \end{align*}
  Using Young's inequality and choosing $\sigma$ sufficiently small, we thus
  obtain that 
  \begin{align*}
    \frac{d}{dt}\langle  \beta, A_1 \beta \rangle \leq \delta \langle  \beta, \dot{A}_1 \beta \rangle  
    + C\sigma^{-1}t^{-2(1-\mu/2)} |\omega_0|_{y=a,b}|^2.
  \end{align*}
  The first part of Theorem \ref{thm:Stabilityofbeta} then follows by integrating this inequality and using a
  Gronwall-type estimate to control the growth up to time $T$.
\end{proof}

Additionally, we make use of the following estimates for boundary evaluations of
$L_t\Gamma, W$ and $\Gamma$, which are obtained as an application of the results
of Appendix \ref{sec:auxil-funct-bound}.
\begin{lem}
\label{lem:bdGamma}
Let $g,h,k$ satisfy the assumptions of the second part of Theorem
\ref{thm:StabilityofGamma}.
Then,
\begin{align*}
  (\p_yL_t\Gamma)(a)=c_1 \langle \Gamma, e^{ikt(y-a)}\tilde{u}_1 \rangle, \\
  (\p_yL_t\Gamma)(b)=c_2 \langle \Gamma, e^{ikt(y-b)}\tilde{u}_2 \rangle, 
\end{align*}
and the following estimates hold:
\begin{align*}
  |\langle \Gamma, e^{ikt(y-a)}\tilde{u}_1 \rangle| &\lesssim \frac{C_\mu}{kt}\sqrt{\sum_n |(\p_y\Gamma)|_{n}^2 <n-kt>^{-2\lambda}} + \frac{C}{kt}|\Gamma(a,t)|, \\
  |\langle W, e^{ikt(y-a)}\tilde{u}_1 \rangle| &\lesssim \frac{C_\mu}{kt}(\|\Gamma\|_{L^2}+\|\beta(t)\|_{L^2})+ \frac{C}{kt}|\omega_0(a,t)|, \\
  |\Gamma(a,t)| &\leq \log(t)(|\omega_0(a)|+\|\omega_0\|_{L^2}).
\end{align*}
\end{lem}

\begin{proof}[Proof of Lemma~\ref{lem:bdGamma}]
  The evaluations of $\p_yL_t\Gamma$ at the boundary are obtained as an
  application of Lemma~\ref{lem:H1calc}.
  The first two estimates follow by integration by parts.
  In order to show the last estimate, we restrict~\eqref{eq:Gamma} to the
  boundary and obtain that 
  \begin{align*}
    |\dt \Gamma(a,t)|\lesssim \frac{1}{kt}(\|\Gamma(t)\|_{L^2}+\|W(t)\|_{L^2}),
  \end{align*}
  where we used that $L_t$ enforces zero Dirichlet data.
  The result hence follows by using Theorem~\ref{thm:L2} and the first part of Theorem
 ~\ref{thm:StabilityofGamma} to control 
  \begin{align*}
    \|\Gamma(t)\|_{L^2}+\|W(t)\|_{L^2}\lesssim |\omega_0(a)|+\|\omega_0\|_{L^2},
  \end{align*}
  and then integrating the inequality.
\end{proof}

\begin{lem}
  \label{lem:commutator}
  Let $W$ be the solution of~\eqref{eq:SLE} with initial datum $\omega_0 \in
  H^1$ and let $\Gamma$ and $\beta$ be as in Lemma~\ref{lem:splitting}.
  Then, for any $\sigma>0$,
  \begin{align*}
  \Re \langle A_1 \p_y\Gamma, \left[\left(\frac{ih'}{k}L_t\cdot + \frac{ih}{k}L_t[\Ell,\p_y]L_t \cdot  \right), \p_y \right] W \rangle \leq \sigma |\langle \p_y\Gamma, \dot{A}_1 \p_y \Gamma \rangle| + C\sigma^{-1} \|W\|_{\tilde{H}^{-1}_t}^2.
  \end{align*}
\end{lem}

\begin{proof}[Proof of Lemma~\ref{lem:commutator}]
  The contribution due to $\frac{ih'}{k}L_t$ can be estimated as in Lemma \ref{lem:com1}.
  In the following we thus focus on the commutator and decompose the commutator into the cases where $\p_y$ falls on $h$,
  \begin{align*}
    \frac{ih''}{k}L_tW + \frac{ih'}{k}L_t[\Ell,\p_y]L_tW,
  \end{align*}
  the terms solving an elliptic equation with vanishing Dirichlet data,
  \begin{align*}
    \frac{ih'}{k}L_t[\Ell,\p_y]L_t W + \frac{ih}{k}L_t[\Ell,\p_y]L_t[\Ell,\p_y]L_tW,
  \end{align*}
  and the homogeneous corrections, 
  \begin{align*}
    \frac{ih'}{k} ((\p_yL_tW)(a,t) e^{ikt(y-a)}u_1+(\p_yL_tW)(b,t) e^{ikt(y-b)}u_2), \\
    \frac{ih}{k} ((\p_yL_t[\Ell,p_y]L_tW)(a,t) e^{ikt(y-a)}u_1+(\p_yL_t[\Ell,\p_y]L_tW)(b,t) e^{ikt(y-b)}u_2),
  \end{align*}

  In the first and second case, we use Lemmata~\ref{lem:compCC}
  and~\ref{lem:EllipticContribution} to estimate by
  \begin{align*}
    \|\p_y\Gamma\|_{\tilde{H}^{-1}_t} \|W\|_{\tilde{H}^{-1}_t},
  \end{align*}
  which is of the desired form by Young's inequality.
\\

  It hence only remains to consider the homogeneous corrections.
Here, we estimate 
\begin{align*}
  \Re \Big\langle  A_1(t)\p_y \Gamma &, \frac{ih'}{k} ((\p_yL_tW)(a,t) e^{ikt(y-a)}u_1+(\p_yL_tW)(b,t) e^{ikt(y-b)}u_2) \\ &+ \frac{ih}{k} ((\p_yL_t[\Ell,p_y]L_tW)(a,t) e^{ikt(y-a)}u_1+(\p_yL_t[\Ell,\p_y]L_tW)(b,t) e^{ikt(y-b)}u_2) \\
&\leq C_\mu \sqrt{\sum_n |\p_y\Gamma_n|^{2}<n-kt>^{-2\mu}} \Big(|\p_yL_tW)(a,t)|+|\p_yL_tW)(b,t)| \\
& \quad +|\p_yL_t[\Ell,p_y]L_tW)(a,t)| + |\p_yL_t[\Ell,p_y]L_tW)(b,t)|\Big). 
\end{align*}
We further recall from Section~\ref{sec:auxil-funct-bound} that boundary
evaluations can be obtained by testing with suitable homogeneous solution to the
adjoint problem.
Hence,
\begin{align*}
  |\p_yL_tW)(a,t)|+|\p_yL_tW)(b,t) \lesssim t^{-1} \|\p_yW\|_{L^2}\leq t^{-1}\|\omega_0\|_{H^1}, \\
  |\p_yL_t[\Ell,p_y]L_tW)(a,t)| + |\p_yL_t[\Ell,p_y]L_tW)(b,t)| \lesssim t^{-1} \|[\Ell,\p_y]L_tW\|_{H^1}.
\end{align*}
We can thus conclude the proof, if we can show that 
\begin{align*}
  \|[\Ell,\p_y]L_tW\|_{H^1} \lesssim \|W\|_{H^1}.
\end{align*}
Expressing $\p_y [\Ell,\p_y]L_tW= [\Ell,\p_y]L_t \p_yW +[[\Ell,\p_y]L_t,\p_y]W$,
this estimate follows from elliptic regularity theory for
$[\Ell,\p_y]L_t|_{t=0}$ and using that multiplication by $e^{ikty}$ is an isometry.
\end{proof}

Building on these results, we can now complete the proof of Theorem
\ref{thm:StabilityofGamma}.
\begin{proof}[Proof of Theorem~\ref{thm:StabilityofGamma}, part 2]
  Following a similar strategy as in the previous part, we consider 
  \begin{align*}
    I_{2}(t):= \langle  \p_y \Gamma, A_1(t)\p_y \Gamma \rangle + C_1 \langle  \Gamma, A_1(t)\Gamma \rangle + C_2 \langle \beta, A_1(t)\beta \rangle + C_3 \langle W, A(t)W\rangle,
  \end{align*}
  where $C_1,C_2,C_3>0$ are to be chosen later.

  Using the preceding results and strategy, it suffices to study 
  \begin{align*}
    \Re \langle \p_t \p_y \Gamma, A_1 \p_y \Gamma \rangle.
  \end{align*}
  Following the same strategy as in the previous part of the proof and using
  Lemma~\ref{lem:bdGamma}, we estimate 
  \begin{align*}
    & \quad \Re \langle  A_1 \p_y\Gamma,   \frac{ih}{k}L_t \p_y \Gamma + \frac{ih'}{k}L_t \Gamma + \frac{ih}{k}L_t[\Ell,\p_y]L_t\Gamma 
\\ & \quad- (\p_yL_t\Gamma)(a) e^{ikt(y-a)}u_1 - (\p_yL_t\Gamma)(b) e^{ikt(y-b)}u_1 \rangle \\
&\leq (C+\sigma+ct^{-(1-\mu)}\log(t)) \|\p_y\Gamma\|_{\tilde{H}^{-1}_t}^2+ \sigma^{-1}|\langle \Gamma, \dot A \Gamma \rangle|, 
  \end{align*}
which can be absorbed.
\\

Furthermore, applying Lemma~\ref{lem:modweight1}, we can control  
\begin{align*}
& \quad\Re \Big\langle  A_1 \p_y\Gamma ,
 \p_y\Big( \frac{h}{k^2t} \langle \Gamma,  e^{ikt(y-a)}\tilde{u}_1 \rangle e^{ikt(y-a)} u_1 - \frac{h}{k^2t} \langle \Gamma,  e^{ikt(y-b)}\tilde{u}_2 \rangle e^{ikt(y-b)} u_2 \\
 & \quad -c_1 \frac{h}{k^2t} \langle W, e^{ikt(y-a)}\p_y\tilde{u}_1 \rangle  e^{ikt(y-a)} u_1  -c_2 \frac{h}{k^2t} \langle W, e^{ikt(y-b)}\p_y\tilde{u}_2 \rangle e^{ikt(y-b)} u_2 \Big) \Big \rangle \\
&\leq C(\mu,g,h,k) \left( \sum_n |\Gamma_n|^2<n-kt>^{-2\mu} \right)^{1/2}
\Big( \left|
   \langle \Gamma,  e^{ikt(y-a)}\tilde{u}_1 \rangle \right| \\
& \quad +\left|\langle \Gamma,  e^{ikt(y-b)}\tilde{u}_2 \rangle\right| +\left|\langle W, e^{ikt(y-a)}\p_y\tilde{u}_1 \rangle\right| 
+\left|\langle W, e^{ikt(y-b)}\p_y\tilde{u}_2 \rangle\right| \Big).
\end{align*}
Applying the estimates of Lemma~\ref{lem:bdGamma} and using Young's inequality,
these contributions can hence again be partially absorbed provided $\sigma$ is
sufficiently small and $T>0$ is sufficiently large.
The remaining non-absorbed terms can be estimated by 
\begin{align*}
  t^{-2(1-\mu/2)} (|\Gamma(a,t)| + |\Gamma(b,t)| + \|\p_y W(t)\|_{L^2} + |\omega_0(a)| + |\omega_0(b)|) \lesssim t^{-2(1-\mu/2)} \|\omega_0\|_{H^1},
\end{align*}
where we used Theorem~\ref{thm:L2}, the first part of Theorem
\ref{thm:StabilityofGamma} and the Sobolev embedding.
\\

It remains to estimate 
\begin{align*}
  \Re \left\langle A_1(t)\p_y\Gamma, \p_y \left( \frac{ih'}{k}L_tW + \frac{ih}{k}L_t[\Ell,\p_y]L_tW  \right) \right\rangle .
\end{align*}
Recalling the definition of $\Gamma$ and $\beta$, we express the right function as 
\begin{align*}
 \left(\frac{ih'}{k}L_t\cdot + \frac{ih}{k}L_t[\Ell,\p_y]L_t \cdot  \right) (\Gamma+\beta) \\
+ \left[\left(\frac{ih'}{k}L_t\cdot + \frac{ih}{k}L_t[\Ell,\p_y]L_t \cdot  \right), \p_y \right] W.
\end{align*}
We then estimate 
\begin{align*}
  \Re \left\langle A_1(t)\p_y\Gamma, \left(\frac{ih'}{k}L_t\cdot + \frac{ih}{k}L_t[\Ell,\p_y]L_t \cdot  \right) (\Gamma+\beta) \right\rangle \lesssim \|\p_y\Gamma\|_{\tilde{H}^{-1}_t} (\|\Gamma\|_{\tilde{H}^{-1}_t} + \|\beta\|_{\tilde{H}^{-1}_t}).
\end{align*}
Using Young's inequality, the respective terms can then again be controlled,
given a suitable choice of $\sigma$.
Finally, using Lemma~\ref{lem:commutator},
\begin{align*}
 \Re \langle A_1 \p_y\Gamma, \left[\left(\frac{ih'}{k}L_t\cdot + \frac{ih}{k}L_t[\Ell,\p_y]L_t \cdot  \right), \p_y \right] W \rangle \\ \leq \sigma |\langle \p_y\Gamma, \dot{A}_1 \p_y \Gamma \rangle| + C\sigma^{-1} \|W\|_{\tilde{H}^{-1}_t}^2, 
\end{align*}
which can again be absorbed and hence concludes the proof.
\end{proof}

\subsection{Weighted stability of $\p_y\beta$ and boundary blow-up}
\label{sec:boundary-blow-up}

In this section we consider the evolution of $\p_y\beta$.
Since the behavior at both boundary points is similar and separates, we for
simplicity of notation consider the case $\omega_0(a)\neq 0$, $\omega_0(b)=0$.
The general case can then be obtained by switching $a$ and $b$ and using the
linearity of the equation.
The function $\beta$ then satisfies~\eqref{eq:beta}:
\begin{align}
  \label{eq:31}
  \begin{split}
  \dt \beta -\frac{ih}{k}L_{t}\beta  - \frac{h}{k^2t}\langle \beta ,e^{ikt(y-a)}u \rangle e^{ikt(y-a)}u &= \omega_{0}(a)\frac{h}{k^2t}e^{ikt(y-a)}u, \\
  \beta|_{t=0}&=0.
  \end{split}
\end{align}
We note that, if $\omega_{0}|_{y=a,b}=0$, then $\beta$ identically vanishes.

We recall that by Proposition
\ref{prop:betaL2} under suitable assumptions on $h,g$ and $k$, $\beta$ is stable
in $L^2$.
However, stability in $H^{1}$ or, indeed in $H^{s},s>1/2$, does not hold due to the asymptotic formation of singularities at the boundary.
\begin{lem}[Boundary blow-up]
  \label{lem:bdblowup}
  Suppose that for some $s>0$, 
  \begin{align*}
\sup_{t>0} \|\beta(t)\|_{H^{s}} = C <\infty.
  \end{align*}
  Then $\beta(a,t)$ satisfies
  \begin{align*}
    |\beta(a,t)-h(a)\omega_{0}(a)k^{-2} \log(t)| \leq C_{s}C,
  \end{align*}
  as $t \rightarrow \infty$
  In particular, if $\omega_{0}(a)\neq 0$, then 
  \begin{align*}
    \sup_{t}\|\beta(t)\|_{C^{0}}=\infty.
  \end{align*}
  Hence, by the Sobolev embedding, in that case,
  \begin{align*}
    \sup_{t} \|\beta(t)\|_{H^{s}} \geq \sup_{t} \log(t) = \infty,
  \end{align*}
  for any $s>\frac{1}{2}$.
\end{lem}

\begin{proof}
  Restricting the evolution by
~\eqref{eq:31} to the boundary, we obtain 
  \begin{align}
    \label{eq:32}
    \dt \beta(a) + \frac{h(a)}{k^2t} \langle \beta, e^{ikt(y-a)}u \rangle &= \frac{h(a)\omega_{0}(a)}{k^2t}.
  \end{align}
  Let $s>0$ and without loss of generality $s<1/2$, then by direct computation
  \begin{align*}
   |\langle \beta, e^{ikt(y-a)}u \rangle | \lesssim C t^{-s}\|\beta\|_{H^s}.
  \end{align*}
  Hence, $\beta(a,t)$ satisfies 
  \begin{align*}
    \dt \beta(a) - \frac{\omega_{0}(a)h(a)}{k^2} \dt \log(t) = t^{-1} \mathcal{O}(t^{-s}) \in L^{1}_{t}.
  \end{align*}
  The result hence follows by integrating in time.
\end{proof}

Letting $s=1$ in the preceding Lemma, we in particular note that in general
$H^1$ stability of $\beta$ fails.
Following a similar approach as in~\cite{Zill5}, one can further show that $s=1/2$ is
indeed critical in the sense that stability holds for $H^{s},s<1/2$.
As this is however not sufficient for optimal decay rates in the damping
estimate of Section~\ref{sec:damping-mixing}, in the following we prove weighted
$H^1$ stability as formulated in Theorem~\ref{thm:Stabilityofbeta}.
Here, we use a different method of proof based Duhamel's formula, the details of
which can be found in Appendix~\ref{sec:Duhamel}.

\subsubsection{Splitting $\p_y\beta$}
\label{sec:weight-h1-stab}

We recall that $\beta$ solves
\begin{align}
  \begin{split}
  \dt \beta- \frac{ih}{k}L_{t}\beta  - \frac{h}{k^2t} \langle \beta, e^{ikt(y-a)}u \rangle e^{ikt(y-a)} u &= \omega_{0}(a) \frac{h}{k^2t}e^{ikt(y-a)}u, \\
  \beta|_{t=0}&=0.
  \end{split}
\end{align}
Applying one $y$ derivative to this equation, we obtain 
\begin{align}
  \label{eq:pybeta}
  \begin{split}
  & \quad \dt \p_{y}\beta - \frac{ih}{k}L_{t}\p_{y}\beta + \frac{h}{k^2t} \langle \p_{y} \beta, e^{ikt(y-a)}u \rangle e^{ikt(y-a)} u
\\ &=  [\frac{ih}{k}L_{t},\p_{y}]\beta + \frac{h}{k^2t} \langle \beta, e^{ikt(y-a)} \p_{y} u \rangle e^{ikt(y-a)} u + \frac{h}{k^2t} \beta(a,t) e^{ikt(y-a)} u \\ &\quad - \frac{h}{k^2t} \langle \beta, e^{ikt(y-a)}u \rangle e^{ikt(y-a)} u 
- \omega_{0}(a) \frac{h'}{k^2t}e^{ikt(y-a)}u \\ & \quad+\omega_{0}(a) \frac{h}{k^2t}e^{ikt(y-a)}\p_{y}u   + \frac{i\omega_{0}(a)h}{k} e^{ikt(y-a)}u, 
  \end{split}
\end{align}
where we used that
\begin{align*}
  \frac{h}{k^2t} \langle \beta, e^{ikt(y-a)}u \rangle \p_y(e^{ikt(y-a)})u= \frac{h}{k^2t} (\langle \p_y(\beta u), e^{ikt(y-a)} \rangle - \beta u e^{ikt(y-a)}|_{y=a}^b) e^{ikt(y-a)}u.
\end{align*}

We note that most terms in~\eqref{eq:pybeta} are very similar to ones in
equation~\eqref{eq:pyGamma} satisfied by $\p_{y}\Gamma$, with the exception of 
\begin{align*}
  \frac{i\omega_{0}(a)h}{k} e^{ikt(y-a)}u,
\end{align*}
which is hence identified as the term driving the blow-up.
Based on this reasoning the following lemma introduces a splitting
of $\p_y \beta$.

\begin{lem}
\label{lem:splitpybeta}
Let $\beta_{I}$ be the solution of 
\begin{align}
  \label{eq:33}
  \begin{split}
  \dt \beta_{I} - \frac{ih}{k}L_{t}\beta_{I} + \frac{1}{ikt} \langle \beta_{I}, e^{ikty}u \rangle e^{ikty} u
&=  [\frac{ih}{k}L_{t},\p_{y}]\beta + \frac{h}{k^2t} \langle \beta, e^{ikty} \p_{y} u \rangle e^{ikty} u \\ & \quad + \frac{h}{k^2t} \beta(a,t) e^{ikty} u - \frac{h}{k^2t} \langle \beta, e^{ikty}u \rangle e^{ikty} u \\
& \quad - \omega_{0}(a) \frac{h}{k^2t}e^{ikty}\p_{y}u +\omega_{0}(a) \frac{h}{k^2t}e^{ikty}\p_{y}u, \\
\beta_{I}|_{t=0}&=0,
  \end{split}
\end{align}
and let $\beta_{II}$ be the solution of
\begin{align*}
  \dt \beta_{II} - \frac{ih}{k}L_{t}\beta_{II} + \frac{h}{k^2t} \langle \beta_{II}, e^{ikty}u \rangle e^{ikty} u &=  \omega_{0}(a) e^{ikty}u, \\
\beta_{V}|_{t=0}&=0.
\end{align*}
Then $\p_y\beta=\beta_{I}+\beta_{II}$.
\end{lem}

Following the same strategy as in Section~\ref{sec:h1gamma}, we obtain $L^{2}$ stability of $\beta_{I}$.
\begin{prop}
\label{prop:beta1}
   Suppose the assumptions of Theorem~\ref{thm:Stabilityofbeta} are satisfied, then 
  \begin{align*}
    \|\beta_{I}(t)\|_{L^{2}} \lesssim |\omega_{0}|_{y=a,b}|.
  \end{align*}
\end{prop}

\begin{proof}
  Following the same strategy as in the proof of Theorem
 ~\ref{thm:StabilityofGamma}, we show that, 
  \begin{align*}
    \frac{d}{dt}\langle \beta_{I}, A_1(t) \beta_{I} \rangle 
\leq \langle \beta_{I}, \dot{A}_1(t) \beta_{I} \rangle + C \|\beta_{I}\|_{\tilde{H}^{-1}_t}^2 \\
+ (Ct^{-(1-\mu)} + \sigma) \sum_n |(\beta_{I})_n|^2<n-kt>^{-2\lambda}t^{-\mu} \\
+ C\sigma^{-1} t^{-2(1-\mu/2)} (|\beta(a,t)|^2 + \|\beta\|_{L^2}^2 + |\omega_0(a)|^2). 
  \end{align*}
  Hence, restricting to $t \geq T>0$ and choosing $\sigma$ sufficiently small, 
  \begin{align*}
    \frac{d}{dt}\langle \beta_{I}, A_1(t) \beta_{I} \rangle \leq C\sigma^{-1} t^{-2(1-\mu/2)} (|\beta(a,t)|^2 + \|\beta\|_{L^2}^2 + |\omega_0(a)|^2) \\ \leq C\sigma^{-1} t^{-2(1-\mu/2)}\log(t)^2 |\omega_0(a)|^2,
  \end{align*}
  where we used Proposition~\ref{prop:betaL2} and that, by equation~\eqref{eq:beta}, 
  \begin{align*}
    |\beta(a,t)|\lesssim \int^t \tau^{-1}\|\beta(\tau)\|_{L^2} d\tau \lesssim \log(t) |\omega_0(a)|.
  \end{align*}
\end{proof}

For later reference, we note that we have thus also proven the following proposition.
\begin{prop}
  Suppose that $g,h,k$ satisfy the assumptions of the second part of Theorem
 ~\ref{thm:StabilityofGamma}.
  Then, for any $\omega_0 \in L^2$, the solution $W$ of~\eqref{eq:SLE} satisfies 
  \begin{align*}
    \|W(t)\|_{H^1} + \|\p_y^2W(t) - \beta_{II}(t)\|_{L^2} \lesssim \|\omega_0\|_{H^2},
  \end{align*}
  where $\beta_{II}$ is given by Lemma \ref{lem:splitpybeta}.
\end{prop}

\begin{proof}
  This result combines Theorems~\ref{thm:L2} and~\ref{thm:StabilityofGamma} and
  Propositions~\ref{prop:betaL2} and~\ref{prop:beta1}.
\end{proof}

\subsubsection{Weighted stability of $\beta_{II}$}
\label{sec:it-works}

In order to complete the proof of Theorem~\ref{thm:Stabilityofbeta}, it only remains to study the stability of 
\begin{align*}
  \dt \beta_{II} - \frac{ih}{k}L_{t}\beta_{II} + \frac{h}{k^2t} \langle \beta_{II}, e^{ikty}u \rangle e^{ikty} u &=  \frac{ih}{k}\omega_{0}(a) e^{ikty}u, \\
\beta_{II}|_{t=0}&=0.
\end{align*}
While it would be possible to study this equation directly, we instead build on
our previous analysis of  
\begin{align}
  \label{eq:34}
   \dt - \frac{ih}{k}L_{t}
\end{align}
and introduce an additional boundary layer $\nu$ (c.f. Theorem~\ref{thm:main1}) solving 
\begin{align}
\label{eq:35}
  (\dt - \frac{ih}{k}L_{t})\nu &= \frac{h}{k}\omega_{0}(a)e^{ikty}, \\
  \nu|_{t=0} &=0,
\end{align}
and also define $\beta_{V}=\beta_{II}-\nu$.
Then $\beta_{V}$ solves
\begin{align}
  \label{eq:36}
  \begin{split}
  \dt \beta_{V} - \frac{ih}{k}L_{t}\beta_{V} + \frac{\langle \beta_{V}, e^{ikty}\tilde{u} \rangle}{ikt} e^{ikty}u &= \frac{\langle \nu, e^{ikty}\tilde{u} \rangle}{ikt} e^{ikty}u. \\
  \beta_{V}|_{t=0}&=0.
  \end{split}
\end{align}

\begin{rem}
  Instead of $\nu$ one might attempt to choose the explicit function
  \begin{align*}
    \int^t \frac{ih}{k}\omega_0(a)e^{ik\tau y} d\tau= \frac{ih}{k}\omega_0(a)\frac{e^{ikty}-1}{iky}=:\chi.
  \end{align*}
  However, we note that part of this function oscillates like $e^{ikty}$ and
  that 
  \begin{align*}
    L_t\chi = e^{ikty}L_0 \frac{h}{k^2y}\omega_0(a) + L_t \frac{h}{k^2y},
  \end{align*}
  where $L_0 \frac{h}{k^2y}\omega_0(a)$ is independent of $t$.
  Hence, even for a constant function $u$
  \begin{align*}
   \langle u, L_t \chi \rangle
  \end{align*}
  would not decay or oscillate rapidly enough to be an integrable perturbation. 
\end{rem}

As the main result of this section we establish the following proposition, which
concludes the proof of Theorem~\ref{thm:Stabilityofbeta}.
\begin{prop}
\label{prop:betaV}
  Suppose the assumptions of Theorem~\ref{thm:Stabilityofbeta} are satisfied.
  Then the functions $\beta_V$ and $\nu$ satisfy
  \begin{align}
    \label{eq:37}
    \|\beta_{V}(t)\|_{L^{2}} \lesssim |\omega_{0}(a)|, \\
    \|(y-a)(y-b)\nu(t)\|_{L^{2}} \lesssim |\omega_{0}(a)|.
  \end{align}
\end{prop}

As the evolution of $\beta_{V}$ depends on $\nu$ via 
\begin{align}
  \label{eq:38}
  \langle \nu, e^{ikt(y-a)}\tilde{u} \rangle
\end{align}
and as our estimates of $\nu$ rely on properties of the solution operator of~\eqref{eq:34} (and hence $W$), we follow a multi-step approach: 
\begin{enumerate}
\item Using Propositions~\ref{prop:beta1} and~\ref{prop:damping}, we show that~\eqref{eq:38} grows at most like $\sqrt{t}$.
\item By direct computation, we show that $\|(y-a)(y-b) \nu\|_{L^{2}}$ grows at most like $\log(t)$.
\item This yields a weaker form of Proposition~\ref{prop:betaV} with an estimate by $\sqrt{t} |\omega_{0}(a)|$.
\item Combining this estimate with the damping result of Section~\ref{sec:damping-mixing}, the estimate of ~\eqref{eq:38} improves to $\log(t)$ and we obtain a uniform bound of $\|(y-a)(y-b) \nu\|_{L^{2}}$.
\item Finally, we establish $L^{2}$ stability of $\beta_{V}$ and thus conclude the proof of Proposition~\ref{prop:betaV}.
\end{enumerate}

\begin{lem}
\label{lem:step1}
Assume that the assumptions of Theorem~\ref{thm:Stabilityofbeta} are satisfied. 
Then~\eqref{eq:38} satisfies 
\begin{align}
  \label{eq:39}
  |\langle e^{ikty}  \tilde u , \nu(t)\rangle | \lesssim \sqrt{t}|\omega_{0}(a)|
\end{align}
as $t \rightarrow \infty$.
\end{lem}

\begin{lem}
\label{lem:step2}
Assume that the assumptions of Theorem~\ref{thm:Stabilityofbeta} are satisfied. 
Then $\nu(t)$ satisfies 
\begin{align}
  \label{eq:40}
  \|(y-a)(y-b)\nu(t)\|_{L^{2}} \lesssim \log(t)|\omega_{0}(a)|
\end{align}
as $t\rightarrow \infty$.
\end{lem}

\begin{lem}
\label{lem:step3}
Assume that the assumptions of Theorem~\ref{thm:Stabilityofbeta} are satisfied.
Then, as $t \rightarrow \infty$, $\beta_{V}$ and $\nu$ satisfy
\begin{align*}
  \|\beta_{V}(t)\|_{L^{2}} &\lesssim \sqrt{t}|\omega_{0}(a)|, \\
  \|(y-a)(y-b)\nu(t)\|_{L^{2}} &\lesssim \log(t)|\omega_{0}(a)|.
\end{align*}
In particular, we conclude that the solution opertor
\begin{align*}
  S(t,0): H^{2}(dy) &\rightarrow H^{2}\left( (y-a)(y-b) dy\right), \\
  \omega_0 &\mapsto W(t),
\end{align*}
  satisfies 
  \begin{align*}
    ||| S(t,0)||| \lesssim \sqrt{t}.
  \end{align*}
\end{lem}

\begin{lem}
\label{lem:step4}
Assume that the assumptions of Theorem~\ref{thm:Stabilityofbeta} are satisfied. 
Then $\nu$ satisfies 
\begin{align}
  \label{eq:41}
  \|(y-a)(y-b)\nu(t)\|_{L^{2}} \lesssim |\omega_{0}(a)|
\end{align}
as $t\rightarrow \infty$.
\end{lem}

\begin{lem}
\label{lem:step5}
Assume that the assumptions of Proposition~\ref{prop:betaV} are satisfied. 
Then $\beta_{V}$ satisfies 
\begin{align}
  \label{eq:42}
  \|\beta_{V}(t)\|_{L^{2}} \lesssim |\omega_{0}(a)|
\end{align}
as $t\rightarrow \infty$.
\end{lem}

In our proof of Lemmata~\ref{lem:step1} to~\ref{lem:step5}, we rely on more
detailed, (semi-explicit) characterization of $\nu(t)$ via Duhamel's formula,
which is established in Appendix~\ref{sec:Duhamel}.

\begin{proof}[Proof of Lemma~\ref{lem:step1}]
We directly compute
\begin{align}
  \label{eq:43}
  \langle e^{ikty}\tilde{u}, \int_{0}^{t} e^{ikty}S(t,\tau) e^{ik\tau y} u d\tau \rangle \\
= \langle \tilde{u}, \int_{0}^{t} e^{ik(t-\tau) y}S(t-\tau,0) u d\tau \rangle .
\end{align}
Next, we integrate 
\begin{align}
  \label{eq:44}
  e^{ik(t-\tau) y} = \p_{\tau} \frac{e^{ik(t-\tau)y}-1}{iky}
\end{align}
by parts in $\tau$.
Here, we obtain a boundary term
\begin{align}
  \label{eq:45}
  \langle \tilde{u},\frac{e^{ikty}-1}{iky} S(t,0) u \rangle
\end{align}
and an integral term 
\begin{align}
  \label{eq:46}
  \langle \tilde{u},\int_{0}^{t}\frac{e^{ik(t-\tau) y}-1}{iky} \p_{\tau} S(t-\tau,0) u d\tau \rangle .
\end{align}
For~\eqref{eq:45} we apply Hölder's inequality and control by 
\begin{align}
  \label{eq:47}
  \|\tilde{u}\|_{L^{\infty}} \|\frac{e^{ikty}-1}{iky}\|_{L^{1}_{y}} \|S(t,0) u\|_{L^{\infty}} \lesssim \log(t) \|u\|_{H^{1}}.
\end{align}
In the integral term we use the damping estimate, Proposition~\ref{prop:damping}, to control by 
\begin{align}
  \label{eq:48}
 \int_{0}^{t} \|\tilde{u}\|_{L^{\infty}}  \|\frac{e^{ik(t-\tau)y}-1}{iky}\|_{L^{2}} \|\p_{\tau} S(t-\tau,0) u\|_{L^{2}} d\tau \\
\lesssim \int_{0}^{t}\sqrt{|t-\tau|} <t-\tau>^{-1} \|S(t-\tau,0) u\|_{H^{1}} d\tau \\
  \lesssim \int <t-\tau>^{-1/2} d\tau \lesssim \sqrt{t}.
\end{align}
\end{proof}

\begin{proof}[Proof of Lemma~\ref{lem:step2}]
Using Lemmata~\ref{lem:tDuhamel} and~\ref{lem:shiftS}, we obtain that 
\begin{align}
  \label{eq:49}
  \nu(t)= \int_{0}^{t} e^{ikt(t-\tau)(y-a)}S(t-\tau,0)u d\tau
\end{align}
Multiplying with $(y-a)$, we use that 
\begin{align*}
  -\p_{\tau} \frac{e^{ik(t-\tau)(y-a)}-1}{ik}= (y-a)e^{ik(t-\tau)(y-a)}
\end{align*}
and hence control 
\begin{align*}
  \|(y-a)\nu(t)\|_{L^{2}} & \leq \| \frac{e^{ik(t-\tau)(y-a)}-1}{ik} S(t-\tau,0)u|_{\tau=0}^{t}\|_{L^{2}} \\
 &\quad + \int_{0}^{t} \|\frac{e^{ik(t-\tau)(y-a)}-1}{ik} \p_{\tau}S(t-\tau,0)u \|_{L^{2}} d\tau \\
&\lesssim |k|^{-1} \|u\|_{L^{2}} + |k|^{-2} \|h\|_{L^{\infty}} \int_{0}^{t} \|L_{t-\tau} S(t-\tau,0)u\|_{L^{2}} \\
& \lesssim |k|^{-1} \|u\|_{L^{2}} + |k|^{-2} \|h\|_{L^{\infty}} \int_{0}^{t} <t-\tau>^{-1} \|S(t-\tau,0)u\|_{H^{1}} d\tau \\
& \lesssim |k|^{-1} \|u\|_{L^{2}} + |k|^{-2} \|h\|_{L^{\infty}} \|u\|_{H^{1}}\log(t),
\end{align*}
where we used Proposition~\ref{prop:damping} and Theorem~\ref{thm:StabilityofGamma}.
\end{proof}

\begin{proof}[Proof of Lemma~\ref{lem:step3}]
Using our Lyapunov functional approach on $\beta_{V}$, we need to estimate 
\begin{align}
\label{eq:50}
  \langle A_{1}\beta_{V}, e^{ikty}u\rangle \frac{\langle \nu, e^{ikty}\tilde{u} \rangle}{ikt}.
\end{align}
By Lemma~\ref{lem:step1}, we control 
\begin{align*}
  |\frac{\langle \nu, e^{ikty}\tilde{u} \rangle}{ikt}|\lesssim t^{-1/2},
\end{align*}
and using Lemma~\ref{lem:A1}, we estimate. 
\begin{align*}
  |\langle A_{1}\beta_{V}, e^{ikty}u\rangle | \leq C_{\lambda} \|(\beta_{V})<n-kt>^{-\lambda}\|_{l^{2}},
\end{align*}
where $0<\lambda<\frac{1}{2}$.

Hence, using Young's inequality, we can control~\eqref{eq:50} by 
\begin{align*}
  \epsilon \|(\beta_{V})_n<n-kt>^{-\lambda}\|_{l^{2}}t^{-1/2} + C(\epsilon,\lambda) t^{-1/2}.
\end{align*}
Here, for $\epsilon$ sufficiently small, the first term can be absorbed by 
\begin{align*}
  \langle \beta_{V}, \dot{A}_{1} \beta_{V} \rangle
\end{align*}
and in summary we obtain 
\begin{align*}
  \dt \langle \beta_{V}, A_{1} \beta_{V} \rangle \leq C(\epsilon,\lambda) t^{-1/2}.
\end{align*}
Integrating this inequality then yields the result.  
\end{proof}
We remark that already in step 3 we could obtain a better growth bound by optimizing in $\lambda$ and the splitting of $t^{-1/2}$ in Young's inequality. However, since $t^{-1/2}\not \in L^{2}$ this would only yield a non-uniform bound and our multi-step proof only requires a better than linear growth bound.  

\begin{proof}[Proof of Lemma~\ref{lem:step4}]
Following the proof of Lemma~\ref{lem:step2} it suffices to show that 
\begin{align}
  \label{eq:51}
  \int_{0}^{t}\|\p_{\tau}S(t-\tau,0)u\|_{L^{2}} d\tau \lesssim 1,
\end{align}
uniformly in $t$.
Using Hölder's inequality and Proposition~\ref{prop:damping}, we estimate 
\begin{align*}
  \|\p_{\tau}S(t-\tau,0)u\|_{L^{2}} &= \|\frac{ih}{k}L_{t-\tau}S(t-\tau,0)u\|_{L^{2}} \leq \|h\|_{L^{\infty}}|k|^{-1} \|L_{t-\tau}S(t-\tau,0)u\|_{L^{2}}  \\
&\leq \|h\|_{L^{\infty}}|k|^{-1}<t-\tau>^{-2}(\|(y-a)(y-b)\p_{y}^{2}S(t-\tau,0)u\|_{L^{2}} \\ & \quad +\|S(t-\tau,0)u\|_{H^{1}}) \leq \|h\|_{L^{\infty}}|k|^{-1}<t-\tau>^{-2} |||S(t-\tau,0)||| \|u\|_{H^{2}},
\end{align*}
the operator norm of $S(t-\tau,0)$ is given by Lemma~\ref{lem:step3}.
Hence, we obtain that
\begin{align*}
  \|\p_{\tau}S(t-\tau,0)u\|_{L^{2}} \lesssim \|h\|_{L^{\infty}}|k|^{-1}\|u\|_{H^{2}}<t-\tau>^{-2} \sqrt{t-\tau},
\end{align*}
which is integrable in $\tau$ and thus concludes the proof.
\end{proof}

\begin{proof}[Proof of Lemma~\ref{lem:step5}]
  We claim that 
  \begin{align}
    \label{eq:52}
    |\langle e^{ikty}  \tilde u , \nu(t)\rangle | \lesssim \log (t)|\omega_{0}(a)|.
  \end{align}
  Following the proof of Lemma~\ref{lem:step3}, this implies that 
  \begin{align}
    \label{eq:53}
    |\langle A_{1}\beta_{V}, e^{ikty}u\rangle \frac{\langle \nu, e^{ikty}\tilde{u} \rangle}{ikt} | \leq C \|(\beta_{V})<n-kt>^{-\lambda}\|_{l^{2}} \frac{\log(t)}{t} \\
\leq \epsilon \|(\beta_{V})<n-kt>^{-\lambda}\|_{l^{2}}^{2} t^{-2\mu} + C(\epsilon) \log(t)^{2}t^{-2(1-\mu)}, 
  \end{align}
  where $C(\epsilon)$ is given by Young's inequality and $0<\mu<1$ is chosen such that $2\lambda+2\mu>1$ and $2(1-\mu)>1$.
  Choosing $\epsilon$ sufficiently small, we thus obtain 
  \begin{align*}
   \dt \langle \beta_{V}, A_{1}\beta_{V} \rangle &\leq \langle \beta_{V}, \dot{A}_{1}\beta_{V} \rangle + \epsilon \|(\beta_{V})<n-kt>^{-\lambda}\|_{l^{2}}^{2} t^{-2\mu} + C(\epsilon) \log(t)^{2}t^{-2(1-\mu)} \\ &\leq C(\epsilon) \log(t)^{2}t^{-2(1-\mu)} \in L^{1}_{t}([1,\infty)).
  \end{align*}
  Integrating this inequality then yields the desired result.
  \\

  It remains to prove the claim~\eqref{eq:52}.
  Here, we estimate 
  \begin{align*}
    |\langle e^{ikty} \tilde{u}, \nu(t) \rangle| \lesssim  \log(t) \|u\|_{H^{1}} + \int_{0}^{t}\|\tilde{u}\|_{L^{\infty}}  \|\frac{e^{ik(t-\tau)y}-1}{iky}\|_{L^{2}} \|\p_{\tau} S(t-\tau,0) u\|_{L^{2}}  d\tau .
  \end{align*}
  Using Lemma~\ref{lem:step3} and Proposition~\ref{prop:damping} we control 
  \begin{align*}
    \|\p_{\tau} S(t-\tau,0) u\|_{L^{2}} \leq <t-\tau>^{-2}|||S(t-\tau,0)||| \|u\|_{H^{2}} \lesssim <t-\tau>^{-3/2} \|u\|_{H^{2}},
  \end{align*}
  and we directly compute that 
  \begin{align*}
    \|\frac{e^{ik(t-\tau)y}-1}{iky}\|_{L^{2}} \lesssim \sqrt{t-\tau}.
  \end{align*}
  Hence, we control
  \begin{align*}
    \int_{0}^{t}\|\tilde{u}\|_{L^{\infty}}  \|\frac{e^{ik(t-\tau)y}-1}{iky}\|_{L^{2}} \|\p_{\tau} S(t-\tau,0) u\|_{L^{2}}  d\tau \\
 \lesssim \|\tilde{u}\|_{L^{\infty}} \|u\|_{H^{2}}  \int_{0}^{t} <t-\tau>^{-1} d\tau \leq \|\tilde{u}\|_{L^{\infty}} \|u\|_{H^{2}} \log(t),
  \end{align*}
  which proves the claim.
\end{proof}
 
\appendix
\section{Auxiliary functions and boundary evaluations}
\label{sec:auxil-funct-bound}
In this section we introduce several
auxiliary functions, which can be used to compute boundary evaluations
of of derivatives of $L_tW$ and related quantities.

\begin{lem}
  \label{lem:homfunctions}
  Let $u_{1},u_{2}$ be solutions of
  \begin{align}
    \label{eq:54}
    \begin{split}
    \Ell u&=0, \\
    z &\in [a,b],
    \end{split}
  \end{align}
  with boundary values
  \begin{align}
    \label{eq:55}
    \begin{pmatrix}
      u_{1}(a) & u_{2}(a) \\ u_{1}(b) & u_{2}(b)
    \end{pmatrix}
                                        =
                                        \begin{pmatrix}
                                          1 & 0 \\ 0 & 1
                                        \end{pmatrix}.
  \end{align}
  Let further $\tilde{u}_{1},\tilde{u}_{2}$
  be solutions to the adjoint problem
  \begin{align}
    \label{eq:56}
    \begin{split}
   \Ell^{*}\tilde{u}:= \left(\left((\frac{\p_{y}}{k}-it)g(y)\right)^{2}- (\frac{\p_{y}}{k}-it)\frac{g(y)}{kr(y)}-\frac{1}{r^{2}(y)}\right)\tilde{u}&=0, \\
    y &\in [a,b],
    \end{split}
  \end{align}
  with boundary values
  \begin{align}
    \label{eq:57}
    \begin{pmatrix}
      \tilde{u}_{1}(a) & \tilde{u}_{2}(a) \\
      \tilde{u}_{1}(b) & \tilde{u}_{2}(b)
    \end{pmatrix}
                         =
                         \begin{pmatrix}
                           1 & 0 \\ 0 & 1
                         \end{pmatrix}.
  \end{align}
  Then $u_{1}, u_{2}, \tilde{u}_{1}, \tilde{u}_{2}$ satisfy
  \begin{align}
    \label{eq:58}
    \begin{split}
    u_{1}(t,r,k)&=e^{ikt(y-a)}u_{1}(0,r,k), \\
    u_{2}(t,r,k)&=e^{ikt(y-b)}u_{2}(0,r,k), \\
    \tilde{u}_{1}(t,r,k)&=e^{ikt(y-a)}\tilde{u}_{1}(0,r,k), \\
    \tilde{u}_{2}(t,r,k)&=e^{ikt(y-b)}\tilde{u}_{2}(0,r,k).
    \end{split}
  \end{align}
\end{lem}

\begin{proof}[Proof of Lemma~\ref{lem:homfunctions}]
  We note that the operators in equations~\eqref{eq:54} and
 ~\eqref{eq:56} are obtained by conjugating by $e^{iktz}$
  and are complex linear.  The result hence follows by noting that
  multiplication by $e^{ikt(y-a)}$
  or $e^{ikt(y-b)}$
  is compatible with the boundary conditions~\eqref{eq:55} and
 ~\eqref{eq:57}.
\end{proof}

\begin{lem}
  \label{lem:H1calc}
  Let $W$ be a given function and let $\Phi$ be a solution of
  \begin{align}
    \label{eq:59}
    \begin{split}
    \Ell \Phi&=W, \\
    \Phi|_{y=a,b}&=0,
    \end{split}
  \end{align}
  and let $u_{1},u_{2},\tilde{u}_{1}, \tilde{u}_{2}$
  be as in Lemma~\ref{lem:homfunctions}.  Define
  \begin{align}
    \label{eq:60}
    \begin{split}
    H^{(1)}&= \frac{k^{2}}{g^{2}(a)}\langle \Phi, \tilde{u}_{1} \rangle_{L^{2}}u_{1} + 
             \frac{k^{2}}{g^{2}(b)}\langle \Phi, \tilde{u}_{2} \rangle_{L^{2}}u_{2}, \\
    \Phi^{(1)}&= \p_{r}\Phi - H^{(1)}
    \end{split}
  \end{align}
  Then $\Phi$ satisfies
  \begin{align}
    \label{eq:61}
    \begin{split}
    \langle W, \tilde{u}_{1} \rangle_{L^{2}}&= \frac{g^{2}(a)}{k^{2}}\p_{r}\Phi(t,k,a) \\
    \langle W, \tilde{u}_{2} \rangle_{L^{2}}&=\frac{g^{2}(b)}{k^{2}}\p_{r}\Phi(t,k,b),
    \end{split}
  \end{align}
  and $\Phi^{(1)}$ solves
  \begin{align}
    \label{eq:62}
    \begin{split}
    \Ell \Phi^{(1)} &=\p_{y}W + \left[\Ell , \p_{y} \right]\Phi, \\
    \Phi^{(1)}|_{y=a,b}&=0.
    \end{split}
  \end{align}
  The function $H^{(1)}$
  is a solution of~\eqref{eq:54} and is called the (first)
  \emph{homogeneous correction}.
\end{lem}

\begin{proof}
 Testing the equation \eqref{eq:59} with the homogeneous solutions of
 Lemma~\ref{lem:homfunctions}, the results follow by integration by parts and
 direct calculations.
\end{proof}

\begin{lem}
  Let $\Phi,W$
  as in Lemma~\ref{lem:H1calc} and let
  $u_{1},u_{2},\tilde{u}_{1}, \tilde{u}_{2}$
  as in Lemma~\ref{lem:homfunctions}.  Then $\Phi$ satisfies
  \begin{align}
    \label{eq:63}
    \begin{split}
    \frac{g^{2}(a)}{k^{2}}\p_{r}^{2}\Phi(t,k,a)&= -\frac{g(a)g'(a)}{k^{2}} \p_{r}\Phi(t,k,a)- \frac{g(y)}{k^{2}r(y)}\p_{y}\Phi(t,k,a) + W(t,k,a), \\
    \frac{g^{2}(a)}{k^{2}}\p_{r}^{2}\Phi(t,k,b)&= -\frac{g(b)g'(b)}{k^{2}} \p_{r}\Phi(t,k,b)- \frac{g(y)}{k^{2}r(y)}\p_{y}\Phi(t,k,b) + W(t,k,b).
    \end{split}
  \end{align}
  Define
  \begin{align}
    \label{eq:64}
    \begin{split}
    H^{(2)}&= \p_{y}^{2}\Phi(t,k,a)u_{1} + \p_{y}^{2}\Phi(t,k,b)u_{2}, \\
    \Phi^{(2)}&= \p_{y}^{2}\Phi- H^{(2)},
    \end{split}
  \end{align}
  then $\Phi^{(2)}$ satisfies
  \begin{align}
    \label{eq:65}
    \begin{split}
    ((g(y)(\frac{\p_{y}}{k}-it))^{2}+ \frac{g(y)}{kr(y)}(\frac{\p_{y}}{k}-it)-\frac{1}{r^{2}(y)})\Phi^{(2)}\\
    =\p_{y}^{2}W + \left[((g(y)(\frac{\p_{y}}{k}-it))^{2}+ \frac{g(y)}{kr(y)}(\frac{\p_{y}}{k}-it)-\frac{1}{r^{2}(y)}), \p_{y}^{2} \right]\Phi, \\
    \Phi^{(2)}|_{y=a,b}=0.
    \end{split}
  \end{align}
  The function $H^{(2)}$
  is a solution of~\eqref{eq:54} and is called the (second)
  \emph{homogeneous correction}.
\end{lem}
\begin{proof}
  Direct computation.
\end{proof}

\section{Duhamel's formula and shearing}
\label{sec:Duhamel}

\begin{lem}[Time dependent Duhamel]
  \label{lem:tDuhamel}
  Let $(L(t))_{t \in \R}$ be a given family of linear operators and
  denote by $S(t,t')$ the solution operator of
  \begin{align*}
    (\partial_t + \frac{ih}{k}L_t) a=0,
  \end{align*}
  mapping a prescribed $a(t')$ to $a(t)$.  Then for any given function
  $F$ the unique solution of
  \begin{align*}
    (\partial_t + \frac{ih}{k}L_t) u&=F,\\
    u(0)&=u_0,
  \end{align*}
  is given by
  \begin{align*}
    u(t)= S(t,0)u_0 + \int_{0}^{t} S(t,t') F(t') dt'.
  \end{align*}
\end{lem}

\begin{proof}
  Since $S(0,0)=Id$, we observe that the such defined $u(t)$ satisfies
  $u(0)=u_0$.  It remains to show that $u$ satisfies the equation.  We
  directly compute
  \begin{align*}
    (\partial_t+ \frac{ih}{k}L_t) u(t)&= (\partial_t+\frac{ih}{k}L_t)S(t,0)u_0 + \int_{0}^{t} (\partial_t +\frac{ih}{k}L_t) S(t,t') F(t') dt' + S(t,t)F(t)\\
                                      &= 0 + \int_{0}^t 0 S(t,t') F(t') dt' +\text{Id} F(t) = F(t).
  \end{align*}
  Here we used that for any $t'$
  \begin{align*}
    (\partial_t +\frac{ih}{k}L_t) S(t,t')=0.
  \end{align*}
  We stress that
  \begin{align*}
    (\partial_t +\frac{ih}{k}L_{\tilde t}) S(t,t')
  \end{align*}
  does not vanish in general for any $\tilde t \neq t$.
\end{proof}

Applying Lemma~\ref{lem:tDuhamel} to~\eqref{eq:35}, we obtain that  
\begin{align}
  \label{eq:66}
  \nu(t)= \omega_{0}(a) \int_{0}^{t}S(t,\tau) e^{ik\tau y}u d\tau,
\end{align}
where $S(t,\tau)$ is the solution operator corresponding to~\eqref{eq:34}.
Since $L_{t}$ was defined by a conjugation of $L_{0}$ with $e^{ikty}$, we can also conjugate $S(t,\tau)$.

\begin{lem}
\label{lem:shiftS}
  Let $\sigma >0$, then for any $0\leq s \leq \tau \leq t$ the solution operator $S$ satisfies 
  \begin{align*}
    S(t,\tau)e^{ik \sigma y} f = e^{ik \sigma y} S(t-\sigma,\tau-\sigma) f
  \end{align*}
  for any $f \in L^{2}$.
\end{lem}

\begin{proof}
  We note that for any $t$ 
  \begin{align*}
    e^{-ik\sigma y} \Ell e^{ik\sigma y} = \mathcal{E}_{t-\sigma}
  \end{align*}
  and that also 
  \begin{align*}
     e^{ik\sigma y} \langle e^{ik \sigma y}f, e^{ikty}u \rangle e^{ikty}u=  \langle f, e^{ik(t-\sigma)y}u \rangle e^{ik(t-\sigma)y}u.
  \end{align*}
  Hence, conjugating the equation by $e^{ikt\sigma y}$ is equivalent to a shift in time, which yields the desired result.
\end{proof}

\bibliographystyle{plain} \bibliography{citations}

\begin{thebibliography}{10}

\bibitem{Bedrossian2015}
Bedrossian, Germain, and Masmoudi.
\newblock On the stability threshold for the {3D} {C}ouette flow in {S}obolev
  regularity.
\newblock {\em arXiv preprint arXiv:1511.01373}.

\bibitem{bedrossian2015dynamics}
Jacob Bedrossian, Pierre Germain, and Nader Masmoudi.
\newblock Dynamics near the subcritical transition of the {3D} {C}ouette flow
  {I}: Below threshold case.
\newblock {\em arXiv preprint arXiv:1506.03720}, 2015.

\bibitem{bedrossian2013asymptotic}
Jacob Bedrossian and Nader Masmoudi.
\newblock Asymptotic stability for the {C}ouette flow in the {2D} {E}uler
  equations.
\newblock {\em Applied Mathematics Research eXpress}, 2014(1):157--175, 2014.

\bibitem{bedrossian2015inviscid}
Jacob Bedrossian and Nader Masmoudi.
\newblock Inviscid damping and the asymptotic stability of planar shear flows
  in the {2D} euler equations.
\newblock {\em Publications math{\'e}matiques de l'IH{\'E}S}, 122(1):195--300,
  2015.

\bibitem{bedrossian2016sobolev}
Jacob Bedrossian, Vlad Vicol, and Fei Wang.
\newblock The sobolev stability threshold for {2D} shear flows near couette.
\newblock {\em arXiv preprint arXiv:1604.01831}, 2016.

\bibitem{Euler_stability}
Freddy Bouchet and Hidetoshi Morita.
\newblock {Large time behavior and asymptotic stability of the {2D} {E}uler and
  linearized {E}uler equations}.
\newblock {\em Physica D: Nonlinear Phenomena}, 239(12):948--966, 2010.

\bibitem{chossat2012couette}
Pascal Chossat and G{\'e}rard Iooss.
\newblock {\em The {C}ouette-{T}aylor Problem}, volume 102.
\newblock Springer Science \& Business Media, 2012.

\bibitem{Lin-Zeng}
Zhiwu Lin and Chongchun Zeng.
\newblock {Inviscid dynamical structures near {C}ouette Flow}.
\newblock {\em Archive for rational mechanics and analysis}, 200(3):1075--1097,
  2011.

\bibitem{stepin1995nonself}
Stanislav~Anatol'evich Stepin.
\newblock Nonself-adjoint {F}riedrichs model in hydrodynamic stability.
\newblock {\em Functional Analysis and Its Applications}, 29(2):91--101, 1995.

\bibitem{Zhang2015inviscid}
Dongyi Wei, Zhifei Zhang, and Weiren Zhao.
\newblock Linear inviscid damping for a class of monotone shear flow in
  {S}obolev spaces.
\newblock {\em arXiv preprint arXiv:1509.08228}, 2015.

\bibitem{Zill3}
Christian Zillinger.
\newblock Linear inviscid damping for monotone shear flows.
\newblock {\em arXiv preprint arXiv:1410.7341}, 2014.

\bibitem{Zill4}
Christian Zillinger.
\newblock {\em Linear inviscid damping for monotone shear flows, boundary
  effects and sharp {S}obolev regularity}.
\newblock PhD thesis, University of Bonn, 2015.

\bibitem{Zill5}
Christian Zillinger.
\newblock Linear inviscid damping for monotone shear flows in a finite periodic
  channel, boundary effects, blow-up and critical {S}obolev regularity.
\newblock {\em Archive for Rational Mechanics and Analysis}, pages 1--61, 2016.

\end{thebibliography}

\end{document}